\documentclass[11pt, twoside, leqno]{article}
\usepackage{amssymb}
\usepackage{amsmath}
\usepackage{amsthm}

\allowdisplaybreaks

\def\rr{{\mathbb R}}
\def\rn{{{\rr}^n}}

\def\zz{{\mathbb Z}}

\def\nn{{\mathbb N}}

\def\cc{{\mathbb C}}

\def\ca{{\mathcal A}}

\def\cd{{\mathcal D}}

\def\cg{{\mathcal G}}
\def\cp{{\mathcal P}}
\def\cq{{\mathcal Q}}

\def\cs{{\mathcal S}}
\def\cx{{\mathcal X}}

\def\az{\alpha}
\def\supp{{\mathop\mathrm{\,supp\,}}}

\def\diam{{\mathop\mathrm{\,diam\,}}}
\def\loc{{\mathop\mathrm{\,loc\,}}}

\def\aoti{{\mathop\mathrm {\,AOTI\,}}}

\def\noz{{\nonumber}}

\def\fz{\infty}
\def\lz{\lambda}
\def\dz{\delta}

\def\wz{\widetilde}
\def\ez{\epsilon}

\def\ezp{{\ez'}}

\def\bz{\beta}

\def\kz{{\kappa}}
\def\gz{{\gamma}}

\def\vz{\varphi}
\def\tz{\theta}
\def\sz{\sigma}

\def\ls{\lesssim}
\def\gs{\gtrsim}

\def\oz{\overline}

\def\wg{\wedge}

\def\lf{\left}
\def\r{\right}

\def\ocg{{\mathring{\cg}}}

\def\lp{{L^p(\cx)}}

\def\df{{\dot F^s_{p,\,q}(\cx)}}

\def\bint{{\ifinner\rlap{\bf\kern.35em--}
\int\else\rlap{\bf\kern.45em--}\int\fi}\ignorespaces}
\def\dbint{\displaystyle\bint}
\def\bbint{{\ifinner\rlap{\bf\kern.35em--}
\hspace{0.078cm}\int\else\rlap{\bf\kern.45em--}\int\fi}\ignorespaces}

\def\f{{F^s_{p,\,q}(\cx)}}

\def\dmsp{{\dot M^{s,p}(\cx)}}
\def\msp{{M^{s,p}(\cx)}}

\def\nj2{{N_{J_2}}}

\def\qji2{{Q^{J_2}_{i_2}}}

\def\qtn{{Q_\tau^{k,\nu}}}
\def\qto{{Q_\tau^{0,\nu}}}

\def\ytn{y_\tau^{k,\nu}}

\def\nkt{{N(k,\tau)}}
\def\ik{{I_k}}

\def\wdk{{\wz D_k}}

\def\ik{{I_k}}

\def\nz{{N(0,\tau)}}

\def\dsum{\displaystyle\sum}
\def\dint{\displaystyle\int}

\def\dsup{\displaystyle\sup}

\newtheorem{thm}{Theorem}[section]
\newtheorem{lem}{Lemma}[section]

\newtheorem{rem}{Remark}[section]
\newtheorem{cor}{Corollary}[section]
\newtheorem{defn}{Definition}[section]

\numberwithin{equation}{section}

\begin{document}

\arraycolsep=1pt

\title{{\vspace{-5cm}\small\hfill\bf J. Funct. Anal., to appear}\\
\vspace{4cm}\bf\Large\bf A Characterization of Haj\l asz-Sobolev
and Triebel-Lizorkin Spaces via Grand Littlewood-Paley Functions\footnotetext {
\hspace{-0.22cm}
2000 {\it Mathematics Subject Classification}.
Primary 42B35; Secondary 42B30, 46E35, 42B25.\endgraf {\it
Key words and phrases}. Sobolev spaces, Triebel-Lizorkin
space, Calder\'on reproducing formula.
\endgraf Dachun Yang was supported by the National
Natural Science Foundation (Grant No. 10871025) of China. Pekka Koskela and
Yuan Zhou were supported by the Academy of Finland grant 120972}}
\author{Pekka Koskela, Dachun Yang and Yuan Zhou}
\date{ }
\maketitle
\begin{center}
\begin{minipage}{13.5cm}
{\small {\bf Abstract} \quad In this paper, we establish the
equivalence between the Haj\l asz-Sobolev spaces or classical
Triebel-Lizorkin spaces and a class of grand Triebel-Lizorkin spaces
on Euclidean spaces and also on metric spaces that are both doubling and
reverse doubling. In particular, when $p\in(n/(n+1),\fz)$,
we give a new characterization of the Haj\l asz-Sobolev spaces $\dot
M^{1,\,p}(\rn)$ via a grand Littlewood-Paley function.
}
\end{minipage}
\end{center}

\bigskip

\section{Introduction\label{s1}}

\hskip\parindent
Recently, analogs of the theory of first order Sobolev spaces on doubling
metric spaces have been established both based on upper gradients
\cite{hk98, ch99, s00} and on pointwise inequalities \cite{h96}. For surveys
on this see \cite{h03,h07}. These different approaches result in the same
function class if the underlying space supports a suitable Poincar\'e
inequality \cite{kz}.
In this paper we further investigate the spaces introduced by
Haj\l asz \cite{h96} (also see \cite{y03}) that are defined via pointwise
inequalities.


\begin{defn}\label{d1.1} Let $(\cx, d)$ be a metric space equipped with a
regular Borel measure $\mu$
such that all balls defined by $d$ have finite and positive
measures.
Let $p\in (0,\fz)$ and $s\in (0,1]$. The
 homogeneous fractional Haj\l asz-Sobolev space $\dmsp$ is the
set of all measurable functions $f\in L^p_\loc(\cx)$
for which there exists a non-negative
function $g\in\lp$ and a set $E\subset\cx$ of measure zero such that
\begin{equation}\label{1.1}
|f(x)-f(y)|\le [d(x,y)]^s[g(x)+g(y)]
\end{equation}
for all $x,\ y\in\cx\setminus E$. Denote by $\cd(f)$ the class of
all nonnegative Borel functions $g$ satisfying \eqref{1.1}.
Moreover, define $\|f\|_\dmsp\equiv\inf_{g\in \cd(f)}\{\|g\|_\lp\},$
where the infimum is taken over all functions $g$ as above.
\end{defn}

In the Euclidean setting, $\dot M^{1,p}$ coincides with the usual homogeneous
first order Sobolev space $\dot W^{1,\,p},$ \cite{h96}, provided $1<p<\infty.$
For
$p\in(n/(n+1),\,1]$, it was very recently proved \cite{ks} that this pointwise
definition yields the corresponding Hardy-Sobolev space. Consequently,
$\dot M^{1,\,p}(\rn)=\dot
F^1_{p,\,2}(\rn)$, for $n/(n+1)<p<\infty,$ where $\dot
F^1_{p,\,q}(\rn)$ refers to a homogeneous Triebel-Lizorkin space (see
Theorem 5.2.3/1 in \cite{t83} and \cite{ks}). In the fractional order,
$s\in (0,1),$ case it was shown in \cite{y03} that
$\dot M^{s,\,p}(\rn)=\dot F^s_{p,\,\fz}(\rn)$, provided $1<p<\infty.$
Notice the jump in the index $q$ when $s$ crosses 1 and that the result in the
fractional order case does not allow for values of $p$ below 1.

We will next introduce a class of grand Triebel-Lizorkin spaces that allow
us to characterize conveniently the fractional Haj\l asz-Sobolev spaces for
$n/(n+s)<p<\infty.$  The definition is based on grand Littlewood-Paley
functions and
we later extend it to the metric space setting, establishing an analogous
characterization.




Let $\zz_+\equiv\nn\cup\{0\}$. Let $\cs(\rn)$ be the Schwartz space,
namely, the space of rapidly decreasing functions endowed with a
family of seminorms
$\{\|\cdot\|_{\cs_{k,\,\,m}(\rn)}\}_{k,\,m\in\zz_+}$, where for any
$k\in\zz_+$ and $m\in(0,\,\fz)$, we set
$$\|\vz\|_{\cs_{k,\,m}(\rn)}\equiv\sup_{\az\in\zz_+^n,\ |\az|\le
k}\sup_{x\in\rn}(1+|x|)^{m}|\partial^\az \vz(x)|.$$ Here we recall
that for any $\az\equiv(\az_1,\cdots,\az_n)\in\zz_+^n$,
$|\az|=\az_1+\cdots+\az_n$ and $\partial^\az\equiv
(\frac{\partial}{\partial x_1})^{\az_1}\cdots
(\frac{\partial}{\partial x_n})^{\az_n}$. It is known that
$\cs(\rn)$ forms a locally convex topology vector space. Denote by
$\cs'(\rn)$ the dual space of $\cs(\rn)$ endowed with the
weak$^\ast$-topology. Moreover, for each $N\in\zz_+$, denote by
$\cs_N(\rn)$ the space of all functions $f\in \cs(\rn)$ satisfying
that $\int_\rn x^\az f(x)\,dx=0$ for all $\az\in\zz_+^n$ with
$|\az|\le N$. For the convenience, we also write
$\cs_{-1}(\rn)\equiv\cs(\rn)$. For any $\vz\in\cs(\rn)$, $t>0$ and
$x\in\rn$, set $\vz_t(x)\equiv t^{-n}\vz(t^{-1}x)$.

For each $N\in\zz_+\cup\{-1\}$, $m\in(0,\,\fz)$ and $\ell\in\zz_+$,
our class of test functions is
\begin{equation}\label{e1.2}
 \ca^\ell_{N,\,m}\equiv\{\phi\in\cs_N(\rn):\
\|\phi\|_{ \cs_{N+\ell+1,\,m}(\rn)}\le1\}.
\end{equation}

\begin{defn}\label{d1.2}
Let $s\in\rr$, $p\in(0,\,\fz)$ and $q\in(0,\,\fz]$.
Let $\ca$ be a class of test functions as in \eqref{e1.2}.
The homogeneous
 grand Triebel-Lizorkin space $\ca\dot F^s_{p,\,q}(\rn)$ is
 defined as the collection of
all $f\in\cs'(\rn)$ such that when $q\in(0,\,\fz)$,
\begin{equation*}
\|f\|_{\ca\dot F^s_{p,\,q}(\rn)}\equiv\lf\|\lf(\sum_{k\in\zz}2^{ksq}
\sup_{\phi\in\ca} |\phi_{2^{-k}}\ast
f|^q\r)^{1/q}\r\|_{L^p(\rn)}<\fz,
\end{equation*}
and when $q=\fz$,
\begin{equation*}
\|f\|_{\ca\dot F^s_{p,\,\fz}(\rn)}\equiv\lf\| \sup_{k\in\zz}2^{ks}
\sup_{\phi\in\ca} |\phi_{2^{-k}}\ast f|\r\|_{L^p(\rn)}<\fz.
\end{equation*}
\end{defn}

For $\ca\equiv \ca^\ell_{N,\,m}$,  we also write $\ca\dot
F^s_{p,\,q}(\rn)$ as $\ca^\ell_{N,\,m}\dot F^s_{p,\,q}(\rn)$.
Moreover, if $N\in\zz_+$ and $\|f\|_{\ca\dot
F^s_{p,\,q}(\rn)}=0$, then it is easy to see that $f\in\cp_N$,
where $\cp_N$ is the space
of polynomials with degree no more than $N$.
So the quotient space $\ca\dot
F^s_{p,\,q}(\rn)/\cp_N$ is a quasi-Banach space.
As usual, an element $[f]=f+\cp_N\in
\ca\dot F^s_{p,\,q}(\rn)/\cp_N$ with $f\in\ca\dot
F^s_{p,\,q}(\rn)$, is simply referred to by $f.$
By abuse of the notation,
we always write the space $\ca\dot
F^s_{p,\,q}(\rn)/\cp_N$ as $\ca\dot
F^s_{p,\,q}(\rn)$.

The grand Triebel-Lizorkin
spaces are closely connected with Haj\l asz-Sobolev spaces
and (consequently) with the classical Triebel-Lizorkin spaces.

\begin{thm}\label{t1.1}
Let $s\in(0,\,1]$ and $p\in(n/(n+s),\,\fz)$. If
$\ca=\ca^\ell_{0,\,m}$ with $\ell\in\zz_+$ and $m\in(n+1,\,\fz)$,
then $\dot M^{s,\,p}(\rn)=\ca\dot F^s_{p,\,\fz}(\rn)$ with
equivalent norms.
\end{thm}

To prove Theorem \ref{t1.1}, for any $f\in L^p(\rn)$, we introduce a
special $g\in \cd(f)$ via a variant of the grand maximal function;
see \eqref{e2.x1} below. When $s=1$, comparing this with the proof
of Theorem 1 of \cite{ks}, we see that the gradient on $f$ appearing
there is transferred to the vanishing moments of the test functions
and the size conditions of the test functions and their first-order
derivatives (see $\ca$) here. We point out that the choice of the
set $\ca$ is very subtle. This is the key point which allows us to
extend Theorem \ref{t1.1} to certain metric measure spaces.
Moreover, to prove Theorem \ref{t1.1}, an imbedding theorem
established by Haj\l asz \cite{h03} is also employed.

Theorem \ref{t1.1} also has a higher-order version.

\begin{defn}\label{d1.4}
Let $p\in (0,\fz)$ and $s\in(k,\,k+1]$ with $k\in\nn$. The
homogeneous Haj\l asz-Sobolev space $\dot M^{s,\,p}(\rn)$ is defined
to be the set of all measurable functions $f\in L^p_\loc(\rn)$ such
that for all $\az\in\zz_+^n$ with $|\az|=k$, $\partial^\az f\in \dot
M^{s-k,\,p}(\rn)$, and normed by $\|f\|_{\dot M^{s,\,p}(\rn)}
\equiv\sum_{|\az|=k} \|\partial^\az f\|_{\dot M^{s-k,\,p}(\rn)}$.
\end{defn}


\begin{cor}\label{c1.1}
Let $N\in\zz_+$, $s\in(N,\,N+1]$ and $p\in(n/(n+N-s),\,\fz)$. If
$\ca=\ca^\ell_{N,\,m}$ with $\ell\in\zz_+$ and $m\in(n+N+2,\,\fz)$
when $s=N+1$ or $m\in(n+N+1,\,\fz)$ when $s\in(N,\,N+1)$, then $\dot
M^{s,\,p}(\rn)=\ca\dot F^{s}_{p,\,\fz}(\rn)$ with equivalent norms.
\end{cor}

The essential point in the proof of Corollary \ref{c1.1} is to establish a
lifting property for $\ca\dot F^{s}_{p,\,\fz}(\rn)$ via Theorem
\ref{t1.1}. This is done with the aid of auxiliary lemmas (see
Lemmas \ref{l2.4} and \ref{l2.5} below), where in Lemma \ref{l2.4},
we decompose a test function in $\cs_N(\rn)$ into a sum of
test functions in $\cs_{k}(\rn)$ with subtle controls on their
semi-norms for all $-1\le k\le N-1$. The decomposition of a test
function in $\cs_0(\rn)$ into functions in $\cs(\rn)$ already plays
a key role in \cite{ks}. The proof of Corollary \ref{c1.1} also
uses Theorem \ref{t1.2} below.

Now we recall the definition of Triebel-Lizorkin spaces on $\rn$.

\begin{defn}\label{d1.3}
Let $s\in\rr$, $p\in(0,\,\fz)$ and $q\in(0,\,\fz]$. Let
$\vz\in\cs(\rn)$ satisfy that
\begin{equation}\label{e1.3}
\supp\widehat\vz\subset \{\xi\in\rn: 1/2\le|\xi|\le2\}\ {and}\
|\widehat\vz(\xi)|\ge {\rm constant}>0 \ {if}\ 3/5\le|\xi|\le5/3.
\end{equation}
\noindent The homogeneous Triebel-Lizorkin space $\dot
F^s_{p,\,q}(\rn)$ is defined as the collection of all
$f\in\cs'(\rn)$ such that
\begin{equation*}
\|f\|_{\dot F^s_{p,\,q}(\rn)}\equiv\lf\|\lf(\sum_{k\in\zz}2^{ksq}
|\varphi_{2^{-k}}\ast f|^q\r)^{1/q}\r\|_{L^p(\rn)}<\fz
\end{equation*}
with the usual modification made when $q=\fz$.
\end{defn}

Notice that if $\|f\|_{\dot
F^s_{p,\,q}(\rn)}=0$,
then it is easy to see that $f\in\cp\equiv\cup_{N\in\nn}\cp_N$.
So similarly to above, we write an element $[f]=f+\cp$ in
the quotient space $\dot
F^s_{p,\,q}(\rn)/\cp$ with $f\in\dot
F^s_{p,\,q}(\rn)$ as $f$,
and also the space $ \dot
F^s_{p,\,q}(\rn)/\cp$ as $ \dot
F^s_{p,\,q}(\rn)$.

\begin{thm}\label{t1.2}
Let $s\in\rr$, $p\in(0,\,\fz)$,
 $q\in(0,\,\fz]$ and
$J\equiv n/\min\{1,\,p,\,q\}$.
If $\ca=\ca^\ell_{N,\,m}$ with $\ell\in\zz_+$,
\begin{equation}\label{e1.4}
N+1>\max\{s,\,J-n-s\}\ {and}\
m>\max\{J,\,n+N+1\},
\end{equation}
then $\ca\dot F^s_{p,\,q}(\rn)=\dot F^s_{p,\,q}(\rn)$ with equivalent
norms.
\end{thm}

To prove Theorem \ref{t1.2}, we use the Calder\'on reproducing
formula in \cite{Pe76,FJ85} and the boundedness of almost diagonal
operators on the sequence spaces corresponding to the
Triebel-Lizorkin spaces. The almost diagonal operators were
introduced by Frazier and Jawerth \cite{FJ90} and proved to be a
very powerful tool therein (see also \cite{BH06}). It is perhaps worthwhile to
point out that the proof of Theorem \ref{t1.1} does not rely on
Theorem \ref{t1.2}.

Recall that $\dot M^{1,\,p}(\rn)= \dot F^1_{p,\,2}(\rn)$ when
$p\in(n/(n+1),\,\fz)$ by \cite{ks}, $\ca\dot
F^0_{p,\,\fz}(\rn)=L^p(\rn)$ when $p\in(1,\,\fz)$
 and $\ca\dot F^0_{p,\,\fz}(\rn)=H^p(\rn)$ when
$p\in(n/(n+1),\,1]$, where $\ca=\ca^{\ell}_{-1,\,m}$ with $\ell\ge1$
and $m\in(n+1,\,\fz)$, and $H^p(\rn)$ is the classical real Hardy
space (see \cite{s93,g2}). Combining these facts with Theorem
\ref{t1.1}, we have the following result.

\begin{cor}\label{c1.2}
(i) If $p\in(n/(n+1),\,\fz)$ and $\ca=\ca^{\ell}_{0,\,m}$ with
$\ell\in\zz_+$ and $m\in(n+1,\,\fz)$, then $\ca\dot
F^1_{p,\,\fz}(\rn)=\dot M^{1,\,p}(\rn)=\dot F^1_{p,\,2}(\rn)$ with
equivalent norms.

(ii) If $s\in(0,\,1)$, $p\in(n/(n+s),\,\fz)$ and
$\ca=\ca^{\ell}_{0,\,m}$ with $\ell\in\zz_+$ and $m\in(n+1,\,\fz)$,
then $\ca\dot F^s_{p,\,\fz}(\rn)=\dot M^{s,\,p}(\rn)=\dot
F^s_{p,\,\fz}(\rn)$ with equivalent norms.

(iii) Let $\ca\equiv\ca^{\ell}_{-1,\,m}$ with $\ell\ge1$ and
$m\in(n+1,\,\fz)$. If $p\in(n/(n+1),\,1]$, then $\ca\dot
F^0_{p,\,\fz}(\rn)=H^p(\rn)=\dot F^0_{p,\,2}(\rn)$ with equivalent
norms; if $p\in(1,\,\fz)$, then $\ca\dot
F^0_{p,\,\fz}(\rn)=L^p(\rn)=\dot F^0_{p,\,2}(\rn)$ with equivalent
norms.
\end{cor}

Moreover, for $N\in\zz_+$, $s\in(N,\,N+1]$ and
$p\in(n/(n+N-s),\,\fz)$, by Corollary \ref{c1.1}, Theorem \ref{t1.2}
and the lifting property of homogeneous Triebel-Lizorkin spaces, we
have that $\dot M^{s,\,p}(\rn)=\dot F^{s}_{p,\,\fz}(\rn)$ with
equivalent norms when $s\in(N,\,N+1)$, and $\dot
M^{N+1,\,p}(\rn)=\dot F^{N+1}_{p,\,2}(\rn)$ with equivalent norms.

\begin{rem} \label{r1.2}
(i) In a sense, Corollary \ref{c1.2} (i) gives a grand maximal
characterizations of Hardy-Sobolev spaces $\dot H^{1,\,p}(\rn)=\dot
F^1_{p,\,2}(\rn)$ with $p\in(n/(n+1),\,1]$, where  $\dot
H^{1,\,p}(\rn)$ is defined as the space of all $f\in \cs'(\rn)$ such
that $\nabla f\in H^p(\rn)$. We point out the advantage of this
grand maximal characterization is that it only depends on the
first-order derivatives of test functions, which can be replaced by
Lipschitz regularity (see Definition \ref{d1.5}). In fact, our approach
transfers the derivatives on $f$ to vanishing
moments, size conditions and Lipschitz regularity of test
functions. This is a key observation, which allows us to extend this
characterization to certain metric measure spaces without any
differential structure; see Theorems \ref{t1.3} and \ref{t1.4}
below.

(ii) We point out that
 Auchser, Russ and Tchamitchian \cite{art05} characterized
the Hardy-Sobolev space $\dot F^1_{p,\,2}(\rn)$ via a maximal
function which is obtained by transferring the gradient on $f$ to a
size condition on the divergence of the vectors formed by certain
test functions; see Theorem 6 of \cite{art05}. However, this
characterization still depends on the derivatives.

(iii) We also point out that Cho \cite{c99} characterized
Hardy-Sobolev spaces $\dot H^{k,\,p}(\rn)=\dot F^k_{p,\,2}(\rn)$ with $k\in\nn$
via a nontangential maximal function by transferring the derivatives
on the distribution to a fixed specially chosen Schwartz function;
see Theorem I of \cite{c99}.

(iv) We finally remark that a continuous version of the grand
Littlewood-Paley function
$(\sum_{k\in\zz}
\sup_{\phi\in\ca}|\phi_{2^{-k}}\ast f|^2)^{1/2}$
with a different choice of $\ca$ was used by Wilson
\cite{w07} to solve a conjecture of R. Fefferman and
E. M. Stein on the weighted boundedness of the classical
Littlewood-Paley $S$-function.
\end{rem}

Finally, we discuss the metric space setting. Let $(\cx,\,d,\,\mu)$ be a metric
measure space. For any $x\in\cx$ and $r>0$, let $B(x,
r)\equiv\{y\in\cx:\, d(x, y)<r\}$. Recall that $(\cx, d, \mu)$ is
called an RD-space if there exist constants $0<C_1\le1\le C_2$ and
$0<\kz\le n$ such that for all $x\in\cx$, $0<r<2\diam(\cx) $ and
$1\le\lz<2\diam(\cx)/r$,
\begin{equation}\label{e1.5}
C_1\lz^\kz\mu(B(x, r))
\le\mu(B(x, \lz r))\le C_2\lz^n\mu(B(x, r)),
\end{equation}
where and in what follows,
$\diam \cx\equiv\sup_{x,\,y\in \cx}d(x,y)$;
see \cite{hmy2}.

We point out that \eqref{e1.5} implies the doubling property,
 there exists a constant $C_0\in[1,\,\fz)$
such that for all $x\in\cx$ and $r>0$,
$\mu(B(x, 2r))\le C_0\mu(B(x, r))$,
and the reverse doubling property:
there exists a constant $a\in(1,\,\fz)$ such that
for all $x\in\cx$ and $0<r<\diam\cx/a$,
$\mu(B(x,\,ar))\ge 2\mu(B(x,\,r))$.
For more equivalent characterizations of RD-spaces and the fact that each
connected doubling space is an RD-space, see \cite{yz09}.

In what follows, we always assume that $(\cx,\,d,\,\mu)$ is an
RD-space. We also assume that $\mu(\cx)=\fz$ in this section and in
Section \ref{s4}. In the remaining part of this section, let
$\ocg(1,\,2)$, $\ocg(x,\,2^{-k},\,1,\,2)$, $(\cg(1,\,2))'$ and
$(\ocg^\ez_0(\bz,\,\gz))'$ be as in Section \ref{s4}.

\begin{defn}\label{d1.5}
Let $s\in(0,\,1]$, $p\in(0,\,\fz)$ and $q\in(0,\,\fz]$.
Let $\ca:=\{\ca_k(x)\}_{x\in\cx,\,k\in\zz}$ and
$ \ca_k(x)= \{\phi\in\ocg(1,\,2),\,
\|\phi\|_{\ocg(x,\,2^{-k},\,1,\,2)}\le1\}$ for all $x\in\cx$.
The homogeneous
grand Triebel-Lizorkin space $\ca\dot F^s_{p,\,q}(\cx)$
is defined to be the set of all
$f\in(\cg(1,\,2))'$ that satisfy
$$\|f\|_{\ca\dot F^s_{p,\,q}(\cx)}\equiv\lf\|\lf\{\sum^\fz_{k=-\fz}
2^{ksq}\sup_{\phi\in \ca_k(\cdot)}
|\langle f,\,\phi\rangle|^q\r\}^{1/q}\r\|_\lp<\fz$$
with the usual modification made when $q=\fz$.
\end{defn}

Here we also point out that $\|f\|_{\ca\dot F^s_{p,\,q}(\cx)}=0$
implies that $f={\rm constant}$. Similarly to the above, we write the
element $[f]=f+\cc$ in
the quotient space $\ca\dot
F^s_{p,\,q}(\cx)/\cc$ with $f\in\ca\dot
F^s_{p,\,q}(\cx)$ as $f$,
and also the space $ \ca\dot
F^s_{p,\,q}(\cx)/\cc$ as $ \ca\dot
F^s_{p,\,q}(\cx)$.

We have the following analog of Theorem \ref{t1.1}.

\begin{thm}\label{t1.3}
Let $s\in (0,1]$ and $p\in(n/(n+s), \fz)$. Then
$\dmsp=\ca\dot F^s_{p,\,\fz}(\cx)$ with equivalent norms.
\end{thm}

The proof of Theorem \ref{t1.3} uses essentially the same ideas as
those used in the proof of Theorem \ref{t1.1}.
For further characterizations of
$\dot M^{1,\,1}(\cx)$ when $\cx$ is
a doubling Riemannian manifold see \cite{bd}.

We recall the definition of homogeneous Triebel-Lizorkin
spaces $\df$ in \cite{hmy2}.

\begin{defn}\label{d1.6}
Let $\ez\in (0,1)$,
 $s\in(0,\,\ez)$, $p\in(n/(n+\ez),\,\fz)$ and
 $q\in(n/(n+\ez),\,\fz].$
Let $\bz,\ \gz\in(0,\,\ez)$ satisfy
\begin{equation}\label{1.6}
\bz\in(s,\,\ez)\ and\ \gz\in(\max\{s-\kz/p,\ n(1/p-1)_+\},\,\ez).
\end{equation}
Let $\{S_k\}_{k\in\zz}$ be an approximation of the
identity of order $\ez$ with bounded support as in Definition \ref{d4.2}.
For $k\in\zz$, set $D_k\equiv S_k-S_{k-1}$.
The homogeneous grand
Triebel-Lizorkin space $\df$ is defined to be the set of all
$f\in(\ocg^\ez_0(\bz,\gz))'$ that satisfy
$$\|f\|_\df\equiv\lf\|\lf\{\sum^\fz_{k=-\fz}
2^{ksq}|D_k(f)|^q\r\}^{1/q}\r\|_\lp<\fz$$ with the usual
modification made when $q=\fz$.
\end{defn}

As shown in \cite{yz09}, the definition of $\df$ is independent
of the choices of $\ez$, $\bz$, $\gz$ and the approximation of the
identity as in Definition \ref{d4.2}.

\begin{thm}\label{t1.4}
Let all the assumptions be as in Definition \ref{d1.6}.
Then $\dot F^s_{p,\,q}(\cx)=\ca\dot F^s_{p,\,q}(\cx)$
with equivalent norms.
\end{thm}

To prove Theorem \ref{t1.4}, we employ the discrete Calder\'on
reproducing formula established in \cite{hmy2}, which was already
proved to be very useful therein.

This paper is organized as follows. In Section \ref{s2}, we present
the proofs of Theorems \ref{t1.1} and \ref{t1.2} and Corollary
\ref{c1.1}. In Section \ref{s3}, we generalize these results to the
inhomogeneous case; see Theorems \ref{t3.1} and \ref{t3.2} and
Corollary \ref{c3.x} below. In Section \ref{s4}, we present the
proofs of Theorems \ref{t1.3} and \ref{t1.4}. Finally, in Section
\ref{s5}, we generalize Theorems \ref{t1.3} and \ref{t1.4} to the
inhomogeneous case; see Theorems \ref{t5.1} and \ref{t5.2} below.

We point out that Theorems \ref{t1.3}, \ref{t1.4}, \ref{t5.1} and
\ref{t5.2} apply in a wide range of settings, for instance, to
Ahlfors $n$-regular metric measure spaces (see \cite{hei}),
$d$-spaces (see \cite{t06}), Lie groups of the polynomial volume
growth (see \cite{v, vsc, nsw, a}), the complete connected
non-compact manifolds with a doubling measure (see
\cite{acdh04,amr08}), compact Carnot-Carath\'eodory (also called
sub-Riemannian) manifolds (see \cite{nsw,ns01,ns04}) and to
boundaries of certain unbounded model domains of polynomial type in
$\cc^N$ appearing in the work of Nagel and Stein (see
\cite{ns06,nsw,ns01,ns04}).

Finally, we state some conventions. Throughout the paper,
we denote by $C$ a positive
constant which is independent
of the main parameters, but which may vary from line to line.
Constants with subscripts, such as $C_0$, do not change
in different occurrences. The symbol $A\ls B$ or $B\gs A$
means that $A\le CB$. If $A\ls B$ and $B\ls A$, we then
write $A\sim B$.
For any $a,\, b\in\rr$, we denote $\min\{a,\, b\}$, $\max\{a,\, b\}$,
and $\max\{a,\, 0\}$
by $a\wg b$, $a\vee b$ and $a_+$, respectively.
If $E$ is a subset of a metric space
$(\cx, d)$, we denote by $\chi_E$ the characteristic
function of $E$.
For any locally integrable function $f$,
we denote by $\bbint_E f\,d\mu$ (or $m_E(f)$) the average
of $f$ on $E$, namely, $\bbint_E f\,d\mu\equiv\frac 1{\mu(E)}\int_E f\,d\mu$.

\section{Proofs of Theorems \ref{t1.1} and \ref{t1.2} and Corollary \ref{c1.1}\label{s2}}

\hskip\parindent To prove Theorem \ref{t1.1}, we need
a Sobolev embedding theorem, which for $s=1$ is
due to Haj\l asz \cite[Theorem 8.7]{h03}, and for $s\in (0,1)$
can be proved by a slight modification of the proof of
\cite[Theorem 8.7]{h03}. We omit the details.

\begin{lem}\label{l2.1}
Let $s\in(0,\,1]$, $p\in(0,\,n/s)$
and $p_\ast=np/(n-sp)$.
Then
there exists a positive constant
$C$ such that for all $u\in \dot M^{s,\,p}(B(x,\,2r))$ and $g\in \cd(u)$,
\begin{eqnarray*}
   \lf(\inf_{c\in\rr}\dbint_{B(x, r)}
    \lf|u(y)-c\r|^{p_\ast}\,dy\r)^{1/p_\ast}
 \le C  r^s\lf(\dbint_{B(x, 2r)}
    [g(y)]^{p}\,dy\r)^{1/p}.
\end{eqnarray*}
\end{lem}

The following result follows from
 Lemma \ref{l2.1}. We omit the details.

\begin{lem}\label{l2.2}
Let $s\in(0,\,1]$, $p\in[n/(n+s),\,n/s)$
and $p_\ast=np/(n-sp)$.
Then for each $u\in \dot M^{s,\,p}(\rn)$,
there exists a constant
$C$ such that
$u-C\in L^{p_\ast}(\rn)$
and $\|u-C\|_{L^{p_\ast}(\rn)}\le \wz C \|u\|_{\dot M^{s,\,p}(\rn)}$,
where  $\wz C$ is a positive
 constant independent of $u$.
\end{lem}

\begin{proof}[Proof of Theorem \ref{t1.1}]
Let $\ca=\ca^\ell_{0,\,m}$ with $\ell\in\zz_+$ and
$m\in(n+1,\,\fz)$. We first prove that if $f\in \ca\dot
F^s_{p,\,\fz}(\rn)$, then $f\in \dot M^{s,\,p}(\rn)$ and
$\|f\|_{\dot M^{s,\,p}(\rn)}\ls \|f\|_{\ca\dot F^s_{p,\,\fz}(\rn)}$.

To see this, we first assume that
 $f$ is a locally integrable function.
Fix $\vz\in\cs(\rn)$ with compact support
 and $\int_\rn \vz(x)\,dx=1$. Notice that
$\vz_{2^{-k}}\ast f(x)\to f(x)$ as $k\to\fz$
for almost all $x\in\rn$.
Then for almost all $x, y\in\rn$, taking $k_0\in\zz$ such that
$2^{-k_0-1}<|x-y|\le 2^{-k_0}$,
we have
\begin{eqnarray*}
        |f(x)-f(y)|&&\le|\vz_{2^{-k_0}}\ast f(x)-
\vz_{2^{-k_0}}\ast f(y)|\\
        &&\quad+\sum_{k\ge k_0}(|\vz_{2^{-k-1}}\ast f(x)-
\vz_{2^{-k}}\ast f(x)|+ |\vz_{2^{-k-1}}\ast f(y)- \vz_{2^{-k}}\ast
f(y)|).
\end{eqnarray*}
Write $\vz_{2^{-k_0}}\ast f(x)-\vz_{2^{-k_0}}\ast f(y)
=(\phi^{(x,\,y)})_{2^{-k_0}}\ast f(x)$ with
$\phi^{(x,\,y)}(z)\equiv\vz(z-2^{k_0}[x-y])-\vz(z)$ and
$\vz_{2^{-k-1}}\ast f(x)- \vz_{2^{-k}}\ast
f(x)=(\vz_{2^{-1}}-\vz)_{2^{-k}}\ast f(x)$. Notice that
$\vz_{2^{-1}}-\vz$ and $\phi^{(x,\,y)}$ are fixed constant multiples
of elements of $\ca^\ell_{0,\,m}$. For all $x\in\rn$, set
\begin{equation}\label{e2.x1}
g(x)\equiv\sup_{k\in\zz}2^{ks}
\sup_{\phi\in\ca^\ell_{0,\,m}}|\phi_{2^{-k}}\ast f(x)|.
\end{equation}
Since $f\in\ca\dot F^s_{p,\,\fz}(\rn)$ and $s\in(0,\,1]$,
we then have $g\in L^p(\rn)$ and
\begin{eqnarray*}
        |f(x)-f(y)|&&\ls|\vz_{2^{-k_0}}\ast f(x)-
\vz_{2^{-k_0}}\ast f(y)|\\
         &&\quad\quad+ \sum_{k\ge k_0}
\sup_{\phi\in\ca^\ell_{0,\,m}}(|\phi_{2^{-k}}\ast f(x)| +
|\phi_{2^{-k}}\ast f(y)|)\\
&&\ls\sum_{k\ge k_0}2^{-ks}[g(x)+g(y)]
  \ls|x-y|^s[g(x)+g(y)].
\end{eqnarray*}
Thus, $f\in \dot M^{s,\,p}(\rn)$ and
$\|f\|_{\dot M^{s,\,p}(\rn)}\ls\|f\|_{\ca\dot F^s_{p,\,\fz}(\rn)}$.

Generally, if $f\in \ca\dot F^s_{p,\,\fz}(\rn)$ is only known to be
an element in $\cs'(\rn)$ at first, then show that we may identify  $f$ with a
locally integrable function $\wz f$ in $\cs'(\rn)$. Indeed, let $\vz$ be as
above. Notice that for all $x\in\rn$, $k\in\zz$ and $i\in\nn$,
\begin{equation}\label{e2.x2}
|\vz_{2^{-k}}\ast f(x)-\vz_{2^{-(k+i)}}\ast f(x)|
\le\sum_{j=0}^{i-1}|\vz_{2^{-k-j}}\ast f(x)- \vz_{2^{-k-j-1}}\ast
f(x)|\ls 2^{-ks}g(x).\end{equation} If $p\in(1,\,\fz)$, then
$\{\vz_{2^{-k}}\ast f-\vz_{2^{-(k+i)}}\ast f\}_{i\in\nn}$ is a
Cauchy sequence in $L^p(\rn)$, which together with the completeness
of $L^p(\rn)$ implies that
 there exists an $f_k\in L^p(\rn)$ such that
$\vz_{2^{-k}}\ast f-\vz_{2^{-(k+i)}}\ast f \to f_k$
 in $L^p(\rn)$ and thus almost everywhere
as $i \to \fz$. Observe that for any $k,\,k'\in\zz$, we have
\begin{eqnarray*}
f_k&&=
\lim_{i\to\fz}[\vz_{2^{-k}}\ast f-
\vz_{2^{-k-i}}\ast f]\\
&&= [\vz_{2^{-k}}\ast f-\vz_{2^{-k'}}\ast f]
+\lim_{i\to\fz}[\vz_{2^{-k'}}\ast f-
\vz_{2^{-k-i}}\ast f]\\
&&= [\vz_{2^{-k}}\ast f-\vz_{2^{-k'}}\ast f]
+f_{k'}
\end{eqnarray*}
in $L^p(\rn)$ and almost everywhere. Set $\wz f\equiv \vz\ast
f-f_0$. Then $\wz f\in L^1_\loc(\rn)$ and $\wz f=\vz_{2^{-k}}\ast
f-f_k$ almost everywhere. Since $\{\vz_{2^{-k}}\ast f\}_{k\in\zz} $
is a sequence of continuous functions that converges to $f$ in
$\cs'(\rn)$ as $k\to\fz$ (see, for example, Lemma 3.8 of
\cite{b03}), then for any $\psi\in\cs(\rn)$, we have
\begin{eqnarray*}
\int_\rn \wz f(x)\psi(x)\,dx&&=
\int_\rn  \vz\ast f(x)\psi(x)\,dx-\lim_{i\to\fz}\int_\rn
[\vz\ast f(x)-\vz_{2^{-i}}\ast f(x)]\psi(x)\,dx\\
&&= \lim_{i\to\fz}\int_\rn
 \vz_{2^{-i}}\ast f(x)\psi(x)\,dx=\langle f,\,\psi\rangle,
\end{eqnarray*}
which implies that $f$ coincides with $\wz f$ in $\cs'(\rn)$. Now we
identify $f$ with the locally integrable function $\wz f$ in
$\cs'(\rn)$. Therefore, by the above proof,
 $\wz f\in\dot M^{s,\,p}(\rn)$ and $\|\wz f\|_{\dot M^{s,\,p}(\rn)}
 \ls\|\wz f\|_{\ca\dot F^s_{p,\,\fz}(\rn)}\sim\|f\|_{\ca\dot F^s_{p,\,\fz}(\rn)}$.
 In this sense, we have that $f\in\dot M^{s,\,p}(\rn)$
 and $$\|f\|_{\dot M^{s,\,p}(\rn)}\ls\|f\|_{\ca\dot F^s_{p,\,\fz}(\rn)}.$$

Now assume that $p\in(n/(n+s), 1]$.
For any $x,\,y\in\rn$, let $k_0\in\zz$ such that
$2^{-k_0-1}<|x-y|\le 2^{-k_0}$.
If $k>k_0$, then
\begin{eqnarray*}
|\vz_{2^{-k}}\ast f(x)-\vz_{2^{-k}}\ast f(y)|
&&\ls |\vz_{2^{-k}}\ast f(x)-
\vz_{2^{-k_0}}\ast f(x)|+ |\vz_{2^{-k}}\ast f(y)-
\vz_{2^{-k_0}}\ast f(y)|\\
&&\quad+
|\vz_{2^{-k_0}}\ast f(x)-\vz_{2^{-k_0}}\ast f(y)|\\
&&\ls |x-y|^s[g(x)+g(y)],
\end{eqnarray*}
where $g$ is as in \eqref{e2.x1}.
If $k\le k_0$, then $2^{k}|x-y|\le1$ and
\begin{eqnarray*}
|\vz_{2^{-k}}\ast f(x)-\vz_{2^{-k}}\ast f(y)|
&&=2^k|x-y||(\phi^{x,\,y})_{2^{-k}}\ast f(x)|\ls
|x-y|^s[g(x)+g(y)],
\end{eqnarray*}
where $g$ is as in \eqref{e2.x1} and for all $z\in\rn$,
$$\phi^{x,\,y}(z)=2^{-k}|x-y|^{-1}[\vz(z)-\vz(z-2^{k}(x-y))].$$
Thus, $\vz_{2^{-k}}\ast f\in \dot M^{s,\,p}(\rn)$ and
$\|\vz_{2^{-k}}\ast f\|_{\dot M^{s,\,p}(\rn)}
\ls\|f\|_{\ca^\ell_{0,\,m}\dot F^s_{p,\,\fz}(\rn)}$ uniformly in
$k\in\zz$. By Lemma \ref{l2.2}, for each $k\in\zz$, there exists a
constant $C_k$ such that $\vz_{2^{-k}}\ast f-C_k\in L^{p_\ast}(\rn)$
with uniform bounded norms. By the weak compactness property of
$L^{p_\ast}(\rn)$, there exists a subsequence which we still denote
by the full sequence such that $\vz_{2^{-k}}\ast f-C_k$ converges weakly in
$L^{p_\ast}(\rn)$ and thus almost everywhere to a certain function
$\widetilde f\in L^{p_\ast}(\rn)$. Moreover, for all $k,\,k'\in\zz$,
since $\vz_{2^{-k}}\ast f- \vz_{2^{-k'}}\ast f\in L^{p}(\rn)$ (see
\eqref{e2.x2}) and $\vz_{2^{-k}}\ast f- \vz_{2^{-k'}}\ast
f+(C_k'-C_k)\in L^{p_\ast}(\rn)$, we know that $C_k'=C_k$. This,
together with the fact that $\vz_{2^{-k}}\ast f\to f \in \cs'(\rn)$
as $k\to\fz$, implies that $f$ coincides with $\wz f +C_0$ in
$\cs'(\rn)$ and hence with $\wz f$ in $\cs'(\rn)/\cc$. Now, we
identify $f$ with the locally integrable function $\wz f$. As in the
case $p\in(1,\,\fz)$, in this case, we also have that $f\in\dot
M^{s,\,p}(\rn)$ and $\| f\|_{\dot M^{s,\,p}(\rn)}\ls\|f\|_{\ca\dot
F^s_{p,\,\fz}(\rn)}.$

Now we show that if $f \in\dot M^{s,\,p}(\rn)$,
then $f\in \ca\dot F^s_{p,\,\fz}(\rn)$
 and $\|f\|_{\ca\dot F^s_{p,\,\fz}(\rn)}
\ls \|f\|_{\dot M^{s,\,p}(\rn)}$.

Let $\phi\in\ca^\ell_{0,\,m}(\rn)$
and $g\in {\mathcal D}(f)$. Then for all
  $k\in\zz$ and $i\ge 0$, by Lemma \ref{l2.1},
  we have that $f\in L^1_\loc(\rn)$ and
\begin{eqnarray*}
    &&\dbint_{B(x,\,2^{-k+i})}
    \lf|f(y)-\dbint_{B(x,\,2^{-k})}f(z)\,dz\r|\,dy\\
        &&\quad\ls\sum_{j=0}^i\dbint_{B(x,\,2^{-k+j})}
    \lf|f(y)-\dbint_{B(x,\,2^{-k+j})}f(z)\,dz\r|\,dy\\
   &&\quad\ls \sum_{j=0}^i 2^{-ks+js}\lf(\dbint_{B(x,\,2^{-k+j})}
    [g(y)]^{n/(n+s)}\,dy\r)^{(n+s)/n}\\
    &&\quad \ls 2^{-ks}2^{i s}\lf[M\lf(g^{n/(n+s)}\r)(x)\r]^{(n+s)/n}.
\end{eqnarray*}
From this, $m>n+1\ge n+s$ and $\int_\rn\phi(x)\,dx=0$, it
follows that for all $k\in\zz$ and $x\in\cx$,
\begin{eqnarray}\label{e2.1}
|\phi_{2^{-k}}\ast f(x)|&&=\lf|\dint_\cx \phi_{2^{-k}}(x-y)
\lf[f(y)-\dbint_{B(x,\, 2^{-k})}f(z)\,dz\r]\,dy\r|\\
 &&\le\sum_{i=0}^\fz 2^{-(m-n)i}
 \dbint_{B(x,\, 2^{-k+i})}\lf|f(y)-\dbint_{B(x, 2^{-k})}
 f(z)\,dz\r|\,dy\noz\\
&&\ls 2^{-ks}\lf[M\lf(g^{n/(n+s)}\r)(x)\r]^{(n+s)/n},\noz
\end{eqnarray}
which together with the $L^{p(n+s)/n}$-boundedness of $M$ implies that
if $p\in (n/(n+s),\fz)$, then
\begin{equation*}
\lf\|\dsup_{k\in\zz}\sup_{\phi\in\ca^\ell_{0,\,m}(\rn)}
2^{ks}|\phi_{2^{-k}}\ast(f)|\r\|_\lp
\ls\lf\|\lf[M\lf(g^{n/(n+s)}\r)\r]^{(n+s)/n}\r\|_\lp
\ls\|g\|_\lp.
\end{equation*}
Moreover, without loss of generality,
we may assume that $M(g^{n/(n+s)})(0)<\fz$.
Then for any $\psi\in\cs(\rn)$, by an argument similar to that
of \eqref{e2.1}, we have that
\begin{eqnarray*}
 \lf|\int_\cx f(x)\psi(x)\,dx\r|&&\le
 \|\psi\|_{L^1(\cx)}\dbint_{B(0,\,1)}|f(z)|\,dz
 +\int_\rn \lf|f(x)-\dbint_{B(0,\,1)}f(z)\,dz\r||\psi(x)|\,dx\\
 &&\ls
 \|\psi\|_{L^1(\cx)}\dbint_{B(0,\,1)}|f(z)|\,dz+
 \|\psi\|_{\cs_{0,\,m}(\rn)}[M(g^{n/(n+s)})(0)]^{(n+s)/n}\\
&&\ls C(f)\|\psi\|_{\cs_{0,\,m}(\rn)},
\end{eqnarray*}
  which implies that $f\in\cs'(\rn)$. Thus,
$f\in \ca\dot F^s_{p,\,\fz}(\rn)$ and
$\|f\|_{\ca\dot F^s_{p,\,\fz}(\rn)}
\ls\|f\|_{\dot M^{s,\,p}(\rn)},$
which completes the proof of Theorem \ref{t1.1}.
\end{proof}

To prove Theorem \ref{t1.2}, we need the following estimate.
\begin{lem}\label{l2.3}
Let $N\in\zz_+\cup\{-1\}$ and $m\in(n+N+1,\,\fz)$. Then there exists
a positive constant $C$ such that for all $x\in\rn$ and
$i,\,j\in\zz$ with $i\ge j$, $\phi\in\cs_N(\rn)$ and
$\psi\in\cs(\rn)$,
$$
|\phi_{2^{-i}}\ast\psi_{2^{-j}}(x)| \le C
\|\phi\|_{\cs_{N+1,\,m}(\rn)}\|\psi\|_{\cs_{0,\,m}(\rn)} 2^{-(i-
j)(N+1)}2^{j n} (1+2^j|x|)^{-m}.$$
\end{lem}

\begin{proof}
Without loss of generality, we may assume that
$\|\phi\|_{\cs_{N+1,\,m}(\rn)}=\|\psi\|_{\cs_{N+1,\,m}(\rn)}=1$.
For simplicty, we only consider the case $N\ge0$. If $j=0$, then by
$\phi\in\cs_N(\rn)$ and the Taylor formula, we have
\begin{eqnarray*}
|\phi_{2^{-i}}\ast\psi(x)|
&&=\lf|\int_\rn\phi_{2^{-i}}(y)\lf\{\psi(x-y)-
\sum_{|\az|\le N}\frac1{\az!}y^\az\partial^\az\psi(x)\r\}\,dy\r|\\
&&\ls\int_{|y|\le (1+|x|)/2} |\phi_{2^{-i}}(y)|\sum_{|\az|=N+1}
|y|^{N+1}|\partial^\az\psi(x-\tz y)|\,dy\\
&&\quad+\int_{|y|> (1+|x|)/2} |\phi_{2^{-i}}(y)\psi(x-y)|\,dy\\
&&\quad+  \int_{|y|> (1+|x|)/2} |\phi_{2^{-i}}(y)|\sum_{|\az|\le
N}|y|^{|\az|} |\partial^\az\psi(x)|\,dy \equiv I_1+I_2+I_3,
\end{eqnarray*}
where $\tz\in[0,\,1]$. If $|y|\le (1+|x|)/2$, then for any
$\tz\in[0,\,1]$, we have that $1+|x|\le 1+|x-\tz y|+|y|$ and hence,
$1+|x|\le 2(1+|x-\tz y|)$. By this and $m\in(n+N+1,\,\fz)$, we
obtain
\begin{eqnarray*}
I_1&&\ls \int_{|y|\le (1+|x|)/2}
\frac{2^{in}|y|^{N+1}}{(1+2^i|y|)^{m}}\frac{1}{(1+|x-\tz y|)^{m}}\,dy\\
&&\ls 2^{-i(N+1)}\frac{1}{(1+|x|)^{m}} \int_\rn
\frac{2^{in}|2^iy|^{N+1}}{(1+2^i|y|)^{m}}\,dy\ls
2^{-i(N+1)}(1+|x|)^{-m}.
\end{eqnarray*}
For $I_2$ and $I_3$, we also have
\begin{eqnarray*}
I_2&&\ls\int_{|y|>(1+|x|)/2}
\frac{2^{in}}{(1+2^i|y|)^{m}}\frac{1}{(1+|x- y|)^{m}}\,dy \ls
2^{-i(N+1)}(1+|x|)^{-m}
\end{eqnarray*}
and
\begin{eqnarray*}
I_3&&\ls\int_{|y|>(1+|x|)/2} \frac{2^{in}}{(1+2^i|y|)^{m}}|y|^{N}
\frac{1}{(1+|x|)^{m}}\,dy\ls 2^{-i(N+1)}(1+|x|)^{-m}.
\end{eqnarray*}
For $j\ne0$, we obtain
\begin{eqnarray*}
|\phi_{2^{-i}}\ast\psi(x)|&&=2^{jn}|\phi_{2^{-(i-j)}}\ast\psi(2^{j}x)|
 \ls
2^{-|i- j|(N+1)}2^{j n} (1+2^{j}|x|)^{-m},
\end{eqnarray*}
which completes the proof of Lemma \ref{l2.3}.
\end{proof}

Now we turn to the proof of Theorem \ref{t1.2}.

\begin{proof}[Proof of Theorem \ref{t1.2}.]

Let $\ca=\ca^\ell_{N,\,m}$ with $\ell\in\zz_+$, $N$ and $m$
satisfying \eqref{e1.4}. Obviously, $\ca\dot F^s_{p,\,q}(\rn)$ is
continuously imbedded into $\dot F^s_{p,\,q}(\rn)$. We now prove
that if $f\in \dot F^s_{p,\,q}(\rn)$, then $f\in \ca\dot
F^s_{p,\,q}(\rn)$ and $\|f\|_{\ca\dot F^s_{p,\,q}(\rn)}\ls
\|f\|_{\dot F^s_{p,\,q}(\rn)}$. This proof is similar to the proof
that the definition of $\dot F^s_{p,\,q}(\rn)$ is independent of
the choice of $\vz$ satisfying \eqref{e1.3}, but a bit more
complicated. In fact, we need to use the boundedness
of almost diagonal operators in sequences spaces. For reader's convenience, we
sketch the argument.

Recall that there exists a function $\psi\in\cs(\rn)$ satisfying the
same conditions as $\vz$ such that
$\sum_{k\in\zz}\widehat{\vz}(2^{-k}\xi)\widehat\psi(2^{-k}\xi)=1$
for all $\xi\in\rn\setminus\{0\}$; see \cite[Lemma (6.9)]{FJW91}.
Then the Calder\'on reproducing formula says that for all
$f\in\cs'(\rn)$, there exist polynomials $P_f$ and
$\{P_i\}_{i\in\nn}$ depending on $f$ such that for all $x\in\rn$,
\begin{equation}\label{e2.x3}
f(x)+P_f(x)=\lim_{i\to-\fz}\lf\{\sum_{j=i}^\fz
 \vz_{2^{-j}}\ast\psi_{2^{-j}}\ast
f (x)+P_i(x)\r\},\end{equation} where the series converges in
$\cs'(\rn)$; see, for example, \cite{Pe76,FJ85}. When $f\in \dot
F^s_{p,\,q}(\rn)$, it is known that the degrees of the polynomials
$\{P_i\}_{i\in\nn}$ here are no more than $\lfloor s-n/p\rfloor$;
see \cite[pp.\,153-155]{FJ90}, and also \cite[p.\,53]{Pe76} and
\cite[pp.\,17-18]{t83}. Recall that $\lfloor\az\rfloor$ for
$\az\in\rr$ denotes the maximal integer no more than $\az$.
Moreover, as shown in  \cite[pp.\,153-155]{FJ90}, $f+P_f$ is
the canonical representative of $f$ in the sense that if
$\vz^{(i)},\,\psi^{(i)}$ satisfy \eqref{e1.3} and
$\sum_{k\in\zz}\widehat{
\vz^{(i)}}(2^{-k}\xi)\widehat{\psi^{(i)}}(2^{-k}\xi)=1$ for all
$\xi\in\rn\setminus\{0\}$ for $i=1,\,2$, then
$P^{(1)}_{f}-P^{(2)}_f$ is a polynomial of degree no more than
$\lfloor s-n/p\rfloor$, where $P^{(i)}_f$ is as in \eqref{e2.x3}
corresponding to $\vz^{(i)},\,\psi^{(i)}$ for $i=1,\,2$. So in this
sense, we identify $f$ with $\wz f\equiv f+P_f$.

Let $\widetilde\vz(x)=\oz{\vz(-x)}$ for all $x\in\rn$. Denote by
$\cq$ the collection of the dyadic cubes on $\rn$. For any dyadic
cube $Q\equiv2^{-j}k+2^{-j}[0,\,1]^n\in\cq$ with certain
$k\in\zz^n$, we set $x_Q\equiv2^{-j}k$, denote by
$\ell(Q)\equiv2^{-j}$ the side length of $Q$ and write
$\vz_Q(x)\equiv2^{j n/2}\vz(2^j x-k)=2^{-jn/2}\vz_{2^{-j}}(x-x_Q)$
for all $x\in\rn$. It is known that for all $x\in\rn$,
\begin{equation}\label{e2.x4}
\vz_{2^{-j}}\ast\psi_{2^{-j}}\ast f (x)=\sum_{\ell(Q)=2^{-j}
}\langle f,\,\wz\vz_Q \rangle \psi_Q(x)
\end{equation}
in $\cs'(\rn)$ and pointwise; see \cite{FJ85,FJW91} and also
\cite[Lemma 2.8]{BH06}. Notice also that $N+1>s$ implies that
$N\ge\lfloor s-n/p\rfloor$. Then for all $f\in \dot
F^s_{p,\,q}(\rn)$, $\phi\in\cs_N(\rn)$ with  $N\ge \lfloor
s-n/p\rfloor$, $i\in\zz$ and $x\in\rn$, by \eqref{e2.x3} and
\eqref{e2.x4}, we have
$$\wz f\ast\phi_{2^{-i}}(x)= \sum_{Q\in\cq}\langle f,\,
\wz\vz_Q \rangle \psi_Q\ast\phi_{2^{-i}}(x) =\sum_{Q\in\cq}t_Q\psi_Q
\ast\phi_{2^{-i}}(x),$$ where $t_Q=\langle f,\,\wz\vz_Q \rangle$,
and by \cite[Theorem 2.2]{FJ90} or \cite[Theorem (6.16)]{FJW91},
\begin{equation}\label{e2.2}
\|f\|_{\dot F^s_{p,\,q}(\rn)}\sim \|\{t_Q\}_{Q\in\cq}\|_{\dot
f^s_{p,\,q}(\rn)}\equiv\lf\|\lf(\sum_{Q\in\cq}[|Q|^{-s/n-1/2}|t_Q|
\chi_Q]^q\r)^{1/q}\r\|_{L^p(\rn)}.
\end{equation}
Moreover, by Lemma \ref{l2.3}, for all $R\in\cq$ with
$\ell(R)=2^{-i}$ and $x\in R$,
 we have
\begin{eqnarray*}
|\wz f\ast\phi_{2^{-i}}(x)|&&\ls\sum_{j\in\zz}
 2^{-|i-j|(N+1)}2^{(i\wedge j)n}2^{-j n/2}
\sum_{\ell(Q)=2^{-j}}\frac{|t_Q|}
 {(1+2^{i\wedge j}|x-x_Q|)^{m}} \\
&&\ls\sum_{j\in\zz}\sum_{\ell(Q)=2^{-j}}
 2^{-|i-j|(n/2+N+1)}
\frac{|t_Q|}
 {(1+2^{i\wedge j}|x_R-x_Q|)^{m}}|R|^{-1/2}.
\end{eqnarray*}
For $R,\,Q\in\cq$ with $\ell(R)=2^{-i}$ and $\ell(Q)=2^{-j}$,
setting
\begin{eqnarray*}
a_{RQ}&&=2^{-|i-j|(n/2+N+1)} (1+2^{i\wedge j}|x_R-x_Q|)^{-m},
\end{eqnarray*}
by \eqref{e1.4}, we have
$$a_{RQ}\le \lf[\frac{\ell(R)}{\ell(Q)}\r]^s
\lf[1+\frac{|x_R-x_Q|}{\max\{\ell(R),\,\ell(Q)\}}\r]^{-J-\ez}
\min\lf\{\lf[\frac{\ell(R)}{\ell(Q)}\r]^{\frac{n+\ez}2},\,
\lf[\frac{\ell(Q)}{\ell(R)}\r]^{J+\frac{\ez-n}2} \r\}$$ for certain
$\ez>0$. Thus $\{a_{RQ}\}_{R,\,Q\in\cq}$ forms an almost diagonal
operator on $\dot f^s_{p,\,q}(\rn)$, which is known to be bounded
on $\dot f^s_{p,\,q}(\rn)$; see \cite[Theorem 3.3]{FJ90} and also
\cite[Theorem (6.20)]{FJW91}. Therefore, by \eqref{e2.2}, we have
\begin{eqnarray*}
\|\wz f\|_{\ca\dot F^s_{p,\,q}(\rn)} &&\ls\lf\|\lf\{\sum_{R\in\cq}
 \lf[|R|^{-s/n-1/2}\lf(\sum_{Q\in\cq}a_{RQ}t_Q\r)
\chi_R\r]^q \r\}^{1/q}\r\|_{L^p(\rn)}\\
&&\ls\lf\| \{t_Q \}_{Q\in\cq}\r\|_{\dot f^s_{p,\,q}(\rn)} \sim
\|f\|_{\dot F^s_{p,\,q}(\rn)},
\end{eqnarray*}
which completes the proof of Theorem \ref{t1.2}.
\end{proof}

To prove Corollary \ref{c1.1}, we need to establish a lifting
property of $\ca\dot F^s_{p,\,q}(\rn)$, which heavily depends on the
following result.

\begin{lem}\label{l2.4}
Let $N\in\zz_+$, $\vz\in\cs_{N}(\rn)$ and $1\le k\le N+1$. Then
there exist functions
$\{\vz_\az\}_{\az\in\zz_+^n,\,|\az|=k}\subset\cs_{N-k}(\rn)$ such
that $\vz=\sum_{|\az|=k} \partial^\az\vz_\az$; moreover, for any
$\ell\in\zz_+$, there exists a positive constant $C$, depending on
$N,\,k,\,\ell$ and $m$,  but not on $\vz$ and $\vz_{\az}$, such that
\begin{equation}\label{e2.3}
\sum_{|\az|=k} \|\vz_\az\|_{\cs_{N+\ell+1,\,m-kn}(\rn)}\le
C\|\vz\|_{\cs_{N+\ell+1,\,m}(\rn)}.
\end{equation}
\end{lem}

\begin{proof}
We begin by proving Lemma \ref{l2.4} for $k=1$. We
point out that when $N=0$, this proof is essentially given by
\cite[Lemma 6]{ks} and \cite[Lemma 3.29]{af03} except for checking the
estimate \eqref{e2.3}. Now assume $N\ge0$. We decompose
 $\vz$ by using the idea appearing in the proof
of \cite[Lemma 6]{ks} and then verify \eqref{e2.3}.

Let $\vz\in \cs_N(\rn)$. We apply induction on $n$. For $n=1$, set
$\psi(x)\equiv\int_{-\fz}^x\vz(y)\,dy$ for all $x\in\rr$. Then
$\vz(x)=\frac{d}{dx}\psi(x)$ for $x\in\rr$. Moreover, for any $0\le
j\le N-1$, by integration by parts, we have $\int_\rr \psi(x)
x^j\,dx=-\frac1{j+1}\int_\rr \vz(x) x^{j+1}\,dx=0$, which means
$\psi\in\cs_{N-1}(\rr)$. Moreover, for all $x\in\rr$, since
$\vz\in\cs_{N}(\rr)$,
$$|\psi(x)|\le\|\vz\|_{\cs_{N+\ell+1,\,m}(\rr)}
\int^{\fz}_{|x|}\frac{1}{(1+|y|)^{m}}\,dy\ls
\|\vz\|_{\cs_{N+\ell+1,\,m}(\rr)}(1+|x|)^{-(m-1)},$$ and for all
$1\le j\le N+\ell+1$,
$$\lf|\frac{d^j}{dx^j}\psi(x)\r|=\lf|\frac{d^{j-1}}{dx^{j-1}}\vz(x)\r|
\le \|\vz\|_{\cs_{N+\ell+1,\,m}(\rr)}(1+|x|)^{-(m-1)}.$$ Thus Lemma
\ref{l2.4} holds for $n=1$.

Suppose that Lemma \ref{l2.4} holds true for a fixed $n\ge1$. Let
$\vz\in\cs_{N}(\rr^{n+1})$.  Without loss of generality, we may
assume that $\|\vz\|_{\cs_{N+\ell+1,\,m}(\rr^{n+1})}=1$.
For any $x\in\rr^{n+1}$, we write $x=(x',\,x_{n+1})$
and define
$h(x')\equiv\int_\rr \vz(x',\,u)\,du,$
where
$x'=(x_1,\,\cdots,\,x_n)\in\rn$.
Then $h\in\cs_{N}(\rn)$. Moreover,
for all $x'\in\rn$ and $\az'\in\zz_+^{n} $
with $|\az'|\le N+\ell+1$, we have
$$|\partial^{\az'}h(x')|\le
\int_\rr\frac{1}{(|1+|x'|+|u|)^m}\,du
\ls\frac{1}{(|1+|x'|)^{m-1}},
$$
which implies that $\|h\|_{\cs_{N+\ell+1, m-1}(\rn)}\ls1$. By
induction hypothesis, we write
$h(x')=\sum_{i=1}^n\frac{\partial}{\partial x_i}h_i(x')$ with
$h_i\in\cs_{N-1}(\rn)$ and
$\|h_i\|_{\cs_{N+\ell+1,\,m-n-1}(\rn)}\ls1$ for $i=1,\,\cdots,\,n$.
Let $a\in\cs(\rr)$ be fixed and $\int_\rr a(u)\,du=1$. For all
$x\in\rr^{n+1}$, set
$\vz_{n+1}(x)\equiv\int_{-\fz}^{x_{n+1}}[\vz(x',\,u)-a(u)h(x')]\,du$
and $\vz_i(x)\equiv a(x_{n+1})h_i(x')$ with $i=1,\,\cdots,\,n$. Then
$\vz_i\in\cs_{N-1}(\rr^{n+1})$ for $i=1,\,\cdots,\,n$. For any
$\ell\le N-1$ and $|\az|\le\ell$ with
$\az=(\az',\,\az_{n+1})\in\zz_+^{n+1}$, by integration by parts
again, we have
 \begin{eqnarray*}
\int_{\rr^{n+1}}
\vz_{n+1}(x)x^\az\,dx
&&=\int_{\rr^{n}}\int_\rr
\int_{-\fz}^{x_{n+1}}\vz(x',\,u) {(x')}^{\az'}x_{n+1}^{\az_{n+1}}
\,du\,dx_{n+1}\,dx'\\
&&=-\frac1{\az_{n+1}+1}\int_{\rr^{n+1}}
 \vz(x){(x')}^{\az'}x_{n+1}^{\az_{n+1}+1}\,dx=0.
 \end{eqnarray*}
For any $\az\in\zz_+^{n+1}$, for $i=1,\,\cdots,\,n$, we have
\begin{eqnarray*}
|\partial ^\az \vz_i(x)|&&\ls
\frac{\|h_i\|_{\cs_{N+\ell+1,\,m-n-1}(\rn)}}
{(1+|x_{n+1}|)^{m-n-1}(1+|x'|)^{m-n-1}}
\ls (1+|x|)^{-(m-n-1)},
\end{eqnarray*}
which implies that $\|\vz_i\|_{\cs_{N+\ell+1,\,m-n-1}(\rr^{n+1})}
\ls1.$ For any $\az\in\zz^{n+1}_+$ with $|\az|\le N+\ell+1$, if
$\az_{n+1}\ne0$, then by $\|h\|_{\cs_{N+\ell+1,\,m-1}(\rn)}\ls1 $,
we have that $|\partial ^\az \vz_{n+1}(x)|\ls
 (1+|x|)^{-(m-n-1)}$ for all $x\in\rn$;
if $\az_{n+1}=0$, then by $\int_\rr [\psi(x',\,u)-a(u)h(x')]\,du=0$,
we have that for all $x\in\rn$,
\begin{eqnarray*}
|\partial ^\az \vz_{n+1}(x)|&&\ls\int_{|x_{n+1}|}^{\fz}
\frac{1}{(1+|x'|+u)^{m}}\,du + \frac{1}{(1+|x'|)^{m-n}}
\int_{|x_{n+1}|}^{\fz}\frac{1}{(1+|u|)^{m}}\,du\\
&&\ls\|\vz\|_{\cs_{N+\ell+1,\,m}
(\rr^{n+1})}(1+|x'|)^{-(m-n-1)}.
\end{eqnarray*}
Thus,  $\|\vz_{n+1}\|_{\cs_{N+\ell+1,\,m-n-1}(\rr^{n+1})} \ls1$,
which completes the proof of Lemma \ref{l2.4}.
\end{proof}

\begin{lem}\label{l2.5}
For any $N$, $\ell\in\zz_+$ and $m,\,m'\in(n+N+2,\,\fz)$,
$\ca^0_{N,\,m'}\dot F^{N+1}_{p,\,\fz}(\rn) =\ca^\ell_{N,\,m}\dot
F^{N+1}_{p,\,\fz}(\rn)$ with equivalent norms.
\end{lem}

\begin{proof}
It suffices to prove that if
$f\in \ca^\ell_{N,\,m}\dot F^{N+1}_{p,\,\fz}(\rn)$,
then $f\in \ca^0_{N,\,m'}\dot F^{N+1}_{p,\,\fz}(\rn)$
and $\|f\|_{\ca^0_{N,\,m'}\dot F^{N+1}_{p,\,\fz}(\rn)}
\ls \|f\|_{\ca^\ell_{N,\,m}\dot F^{N+1}_{p,\,\fz}(\rn)}.$
Without loss of generality, we may assume that $m\ge m'$.
To this end, fix $\psi\in\cs(\rn)$ such that $\int_\rn \psi(x)\,dx =1$.
Obviously, for any $\az\in\zz_+^n$ with $|\az|=N+1$,
if $\bz\in\zz_+^n$, $|\bz|\le N+1$ and $\bz\ne \az$, then
 $\int_\rn \partial^\az\psi(x)x^\bz\,dx=0$;
if $\az=\bz$, then $\int_\rn
\partial^\az\psi(x)x^\az\,dx=(-1)^{N+1}$. For any $\phi\in
\ca^0_{N,\,m'}$, let
\begin{equation}\label{e2.x5}
\oz\phi=\phi-(-1)^{N+1}\sum_{|\az|=N+1}\lf(\int_\rn
\phi(x)x^\az\,dx\r)
\partial^\az\psi.\end{equation}
Then $\oz\phi\in \cs_{N+1}(\rn)$. Moreover, for $|\az|=N+1$, since
$\phi\in\ca^0_{N,\,m'}$ with $m'\in(n+N+2,\,\fz)$, we have
$$\int_\rn |\phi(x)x^\az|\,dx \le\int_\rn
\frac{|x|^{N+1}}{(1+|x|)^{m'}}\,dx\ls1,$$ which implies that
$\oz\phi$ is a fixed constant multiple of an element of
$\ca^0_{N+1,\,m'}$. Notice that $\{\partial^\az\psi\}_{|\az|=N+1}$
are also fixed constant multiples of elements of $\ca^\ell_{N,\,m}$.
Then, by \eqref{e2.x5}, we have
$$\sup_{\phi\in\ca^0_{N,\,m'} }|\phi_{2^{-k}}\ast f(x)|
\ls\sup_{\phi\in\ca^\ell_{N,\,m} }|\phi_{2^{-k}}\ast f(x)|+
\sup_{\phi\in\ca^0_{N+1,\,m'} }|\phi_{2^{-k}}\ast f(x)|,$$ which
implies that $\|f\|_{\ca^0_{N,\,m'}\dot F^{N+1}_{p,\,\fz}(\rn)} \ls
\|f\|_{\ca^\ell_{N,\,m}\dot F^{N+1}_{p,\,\fz}(\rn)}+
\|f\|_{\ca^0_{N+1,\,m'}\dot F^{N+1}_{p,\,\fz}(\rn)}. $ By Theorem
\ref{t1.2} together with $m'\in(n+N+2,\,\fz)$, we have that
$\|f\|_{\ca^0_{N+1,\,m'}\dot F^{N+1}_{p,\,\fz}(\rn)} \sim
\|f\|_{\dot F^{N+1}_{p,\,\fz}(\rn)}\ls \|f\|_{\ca^\ell_{N,\,m}\dot
F^{N+1}_{p,\,\fz}(\rn)},$ which yields that
$\|f\|_{\ca^0_{N,\,m'}\dot F^{N+1}_{p,\,\fz}(\rn)} \ls
\|f\|_{\ca^\ell_{N,\,m}\dot F^{N+1}_{p,\,\fz}(\rn)}$. This finishes
the proof of Lemma \ref{l2.5}.
\end{proof}

\begin{proof}[Proof of Corollary \ref{c1.1}.]
First, let $f\in \dot M^{s,\,p}(\rn)$. Then by Theorem \ref{t1.1},
for any $\ell\in\zz_+$ and $m\in((n+2)N+1,\,\fz)$, and all
$\az\in\zz_+^n$ with $|\az|=N$, $\partial^\az f\in \dot
M^{s-N,\,p}(\rn)= \ca^\ell_{0,\,m-nN}\dot F^{s-N}_{p,\,\fz}(\rn)$.
Moreover, for any $\phi\in\ca^\ell_{N,\,m}$, by Lemma \ref{l2.4},
there exist $\{\phi_\az\}_{|\az|=N}$ and a positive constant $C$
independent of $\phi$ such that
$\{\frac1C\phi_\az\}_{|\az|=N}\subset \ca^{\ell}_{0,\,m-nN}$ and
$\phi=\sum_{|\az|=N}\partial^\az\phi_\az$. This implies that for all
$x\in\rn$,
$$\phi_{2^{-k}}\ast f(x)=\sum_{|\az|=N}(\partial^\az\phi_\az)_{2^{-k}}\ast f(x)
= 2^{kN}(-1)^N\sum_{|\az|=N}(\phi_\az)_{2^{-k}}\ast (\partial^\az
f)(x),$$ and thus,
$$\sup_{\phi\in \ca^{\ell}_{N,\,m}}|\phi_{2^{-k}}\ast
f(x)|\ls 2^{kN}\sup_{|\az|=N}\sup_{\phi\in
\ca^{\ell}_{0,\,m-nN}}|\phi_{2^{-k}}\ast (\partial^\az f)(x)|.$$ From this
and
 $\partial^\az f\in
\ca^{\ell}_{0,\,m-nN}\dot F^{s-N}_{p,\,\fz}(\rn)$ for all
$\az\in\zz_+^n$ with $|\az|=N$ together with Theorem \ref{t1.1}, it
follows that
  $f\in\ca^{\ell}_{N,\,m}\dot F^{s}_{p,\,\fz}(\rn)$
and $$\|f\|_{\ca^{\ell}_{N,\,m}\dot F^{s}_{p,\,\fz}(\rn)}
\ls\sum_{|\az|=N} \|\partial^\az f\|_{\ca^{\ell}_{0,\,m-nN}
\dot F^{s-N}_{p,\,\fz}(\rn)}
\sim\sum_{|\az|=N}\|\partial^\az f\|_{\dot M^{s-N,\,p}(\rn)}
\sim \|f\|_{\dot M^{s,\,p}(\rn)}.$$

On the other hand, let $f\in\ca^{\ell}_{N,\,m}\dot
F^{s}_{p,\,\fz}(\rn)$. Let $\ell\ge N$ and $m\in(n+N+1,\,\fz)$.
Observe that for any $\phi\in\ca^\ell_{0,\,m}$ and $\az\in\zz_+^n$
with $|\az|=N$, $\partial^\az \phi\in \ca^{\ell-N}_{N,\,m}$. Thus
for all $k\in\zz$,
$$\sup_{\phi\in\ca^{\ell}_{N,\,m}}|\phi_{2^{-k}}\ast(\partial^\az f)|
\le
\sup_{\phi\in\ca^{\ell-N}_{0,\,m}}2^{kN}|\phi_{2^{-k}}\ast(f)|,$$
which implies that $\{\partial^\az
f\}_{|\az|=N}\subset\ca^{\ell-N}_{0,\,m} \dot
F^{s-N}_{p,\,\fz}(\rn)=\dot M^{s-N,\,p}(\rn)$ and thus, $f\in
M^{s,\,p}(\rn)$ and
$$
\|f\|_{\dot M^{s,\,p}(\rn)}\sim
\sum_{|\az|=N}\|\partial^\az f\|_{\dot M^{s-N,\,p}(\rn)}\sim
\sum_{|\az|=N} \|\partial^\az f\|_{\ca^{\ell-N}_{0,\,m}
\dot F^{s-N}_{p,\,\fz}(\rn)}
\ls\|f\|_{\ca^{\ell}_{N,\,m}\dot F^s_{p,\,\fz}(\rn)}.$$

Finally, combining the above results with Lemma \ref{l2.5} and
Theorem \ref{t1.2}, for all $\ell\in\zz_+$ and $m\in(n+N+2,\,\fz)$
when $s=N+1$ or $m\in(n+N+1,\,\fz)$ when $s\in(N,\,N+1)$, we have
that $\dot M^{s,\,p}(\rn)= \ca^{\ell}_{N,\,m}\dot
F^s_{p,\,\fz}(\rn)$. This finishes the proof of Corollary
\ref{c1.1}.
\end{proof}

\section{Inhomogeneous versions of Theorems \ref{t1.1}
and \ref{t1.2}\label{s3}}

\hskip\parindent We first recall the definitions of inhomogeneous
Triebel-Lizorkin spaces; see \cite{t83}.

\begin{defn}\label{d3.1}
Let $s\in\rr$, $p\in(0,\,\fz)$ and $q\in(0,\,\fz]$.
Let $\vz\in\cs(\rn)$ satisfy \eqref{e1.3} and
$\Phi\in\cs(\rn)$ be such that $\supp\widehat\Phi\subset B(0,\,2)$ and
$|\Phi(\xi)|\ge {\rm constant}>0$ for all $|\xi|\le 5/3$.
The inhomogeneous
Triebel-Lizorkin space $F^s_{p,\,q}(\rn)$ is defined as the collection of
all $f\in\cs'(\rn)$ such that
\begin{equation*}
\|f\|_{ F^s_{p,\,q}(\rn)}\equiv\|\Phi\ast f\|_{L^p(\rn)}
+\lf\|\lf(\sum_{k=1}^\fz2^{ksq}
|\varphi_{2^{-k}}\ast f|^q\r)^{1/q}\r\|_{L^p(\rn)}<\fz
\end{equation*}
with the usual modification made when $q=\fz$.
\end{defn}

Recall that the local Hardy space $h^p(\rn)$ of Goldberg is just
$F^0_{p,\,2}(\rn)$ (see \cite[Theorem 2.5.8/1]{t83}). A variant of
inhomogeneous Haj\l asz-Sobolev spaces is defined as follows.

\begin{defn}\rm\label{d3.2}
Let $p\in (0,\fz)$ and $s\in (0,1]$. The
 inhomogeneous fractional Haj\l asz-Sobolev space
$ M^{s,\,p}(\rn)$ is the set of all measurable functions $f\in
L^p_\loc(\rn)$ such that $f\in \dot M^{s,\,p}(\rn)$ and $f\in
h^p(\rn)$. Moreover, define $\|f\|_{ M^{s,\,p}(\rn)}\equiv
\|f\|_{h^p(\rn)}+\|f\|_{\dot M^{s,\,p}(\rn)}.$
\end{defn}

Now we introduce the inhomogeneous grand Triebel-Lizorkin spaces.

\begin{defn}\label{d3.3}
Let $s\in\rr$, $p\in(0,\,\fz)$ and $q\in(0,\,\fz]$. Let
$\ca=\ca^\ell_{N,\,m}$ with $\ell\in\zz_+$, $N\in\zz_+\cup\{-1\}$
and $m\in(0,\,\fz)$
 be a class of test functions as in \eqref{e1.2}.
The inhomogeneous
 grand Triebel-Lizorkin space $\ca F^s_{p,\,q}(\rn)$ is
 defined as the collection of
all $f\in\cs'(\rn)$ such that
\begin{eqnarray}\label{e.1}
\|f\|_{\ca F^s_{p,\,q}(\rn)}&&\equiv
\lf\|\sup_{\phi\in\ca^{\ell+1}_{-1,\,m}}
|\phi\ast f|\r\|_{L^p(\rn)}
+
\lf\|\lf(\sum_{k=1}^\fz2^{ksq}
\sup_{\phi\in\ca}
|\phi_{2^{-k}}\ast f|^q\r)^{1/q}\r\|_{L^p(\rn)}<\fz\nonumber
\end{eqnarray}
with the usual modification made when $q=\fz$.
\end{defn}

Moreover, similarly to Theorems \ref{t1.1} and \ref{t1.2} and
Corollary \ref{c1.1}, we have the following results.

\begin{thm}\label{t3.1}
Let all the assumptions be the same as in Theorem \ref{t1.2}.
Then $\ca F^s_{p,\,q}(\rn)= F^s_{p,\,q}(\rn)$ with
equivalent norms.
\end{thm}

The proof of Theorem \ref{t3.1} is similar to that of Theorem
\ref{t1.2}. In fact, since the inhomogeneous Calder\'on reproducing
formula is available (see \cite[p.\,131]{FJ90}), then by using the
argument of Theorem \ref{t1.2} and the estimates in Lemma
\ref{l2.3}, we have Theorem \ref{t3.1}. Here we omit  the details.

\begin{thm}\label{t3.2}
Let $s\in(0,\,1]$, $p\in(n/(n+s),\,\fz)$, $\ell\in\zz_+$ and
$m\in(n+1,\,\fz)$. Then $M^{s,\,p}(\rn)=\ca^\ell_{0,\,m}
F^s_{p,\,\fz}(\rn)$ with equivalent norms.
\end{thm}

\begin{proof} The proof of Theorem \ref{t3.2}
is similar to that of Theorem \ref{t1.1} and much easier. In fact,
if $f\in M^{s,\,p}(\rn)$, then $ f\in \dot M^{s,\,p}(\rn)$ by
Definition \ref{d3.1} and thus $f\in \ca^\ell_{0,\,m} \dot
F^s_{p,\,\fz}(\rn)$. Notice that
$\|\sup_{\phi\in\ca^{\ell}_{-1,\,m}}|\phi\ast f|\|\ls
\|f\|_{h^p(\rn)}$ (see \cite{g}). Then we know that $f\in
\ca^\ell_{0,\,m} F^s_{p,\,\fz}(\rn)$ and $\|f\|_{\ca^\ell_{0,\,m}
F^s_{p,\,\fz}(\rn)}\ls\|f\|_{M^{s,\,p}(\rn)}$.

Conversely, assume that $f\in \ca^\ell_{0,\,m} F^s_{p,\,\fz}(\rn)$.
Obviously, by $h^p(\rn)=F^0_{p,\,2}(\rn)$, we know that $f\in
h^p(\rn)$ and $\|f\|_{h^p(\rn)}\sim \|f\|_{F^0_{p,\,2}(\rn)}\ls
\|f\|_{\ca^\ell_{0,\,m}F^1_{p,\,2}(\rn)}.$ Fix $\vz\in\cs(\rn)$ with
compact support
 and $\int_\rn \vz(x)\,dx=1$. For any $k\in\nn$,
if $|x-y|\le 1$, by an argument similar to that for Theorem
\ref{t1.1},  we then know that
$$|\vz_{2^{-k}}\ast f(x)-\vz_{2^{-k}}\ast f(y)|
\ls |x-y|^s\sup_{k\ge0}2^{ks}\sup_{\phi\in\ca^{\ell}_{0,\,m}}
(|\phi_{2^{-k}}\ast f(x)|+|\phi_{2^{-k}}\ast f(y)|).$$
If $|x-y|>1$, then, obviously,
$$|\vz_{2^{-k}}\ast f(x)-\vz_{2^{-k}}\ast f(y)|\ls
|x-y|^s \sup_{\phi\in\ca^{\ell+1}_{-1,\,m}} (|\phi_{2^{-k}}\ast
f(x)|+|\phi_{2^{-k}}\ast f(y)|).$$ So $\vz_{2^{-k}}\ast f\in \dot
M^{s,\,p}(\rn)$ and $\|\vz_{2^{-k}}\ast f\|_{\dot M^{s,\,p}(\rn)}\ls
\|f\|_{\ca^\ell_{0,\,m}F^s_{p,\,q}(\rn)}.$ Then similarly to the
proof of Theorem \ref{t1.1}, we can prove that $f\in L^p_\loc(\rn)$
and $\|\vz_{2^{-k}}\ast f\|_{\dot M^{s,\,p}(\rn)}\ls
\|f\|_{\ca^\ell_{0,\,m}F^s_{p,\,q}(\rn)}.$ Thus,
$$\|f\|_{M^{s,\,p}(\rn)}\ls \|f\|_{h^p(\rn)}+
\|f\|_{\ca^\ell_{0,\,m}F^s_{p,\,q}(\rn)}\ls
\|f\|_{\ca^\ell_{0,\,m}F^s_{p,\,q}(\rn)},$$ which completes the
proof of Theorem \ref{t1.3}.
\end{proof}

\begin{cor}\label{c3.1}
Let $s\in(0,\,1]$, $p\in(n/(n+s),\,\fz)$,
 $\ell\in\zz_+$ and $m\in(n+1,\,\fz)$.

(i) If $s=1$, then  $\ca^\ell_{0,\,m} F^1_{p,\,\fz}(\rn)
= M^{1,\,p}(\rn)= F^1_{p,\,2}(\rn)$ with equivalent norms.

(ii) If $s\in (0,1)$, then $\ca^\ell_{0,\,m} F^s_{p,\,\fz}(\rn)
= M^{s,\,p}(\rn)= F^s_{p,\,\fz}(\rn)$ with equivalent norms.
\end{cor}

Define the inhomogeneous Haj\l asz-Sobolev spaces $M^{s,\,p}(\rn)$
of higher orders as in Definition \ref{d1.4} by replacing $\dot
M^{s-N,\,p}(\rn)$ with $M^{s-N,\,p}(\rn)$. Then we have the
following inhomogeneous version of Corollary \ref{c1.1}. We omit the
details of its proof.

\begin{cor}\label{c3.x}
Let $N\in\zz_+$, $s\in(N,\,N+1]$ and $p\in(n/(n+s-N),\,\fz)$. If
$\ca=\ca^\ell_{N,\,m}$ with $\ell\in\zz_+$ and $m\in(n+N+2,\,\fz)$
when $s=N+1$ or $m\in(n+N+1,\,\fz)$ when $s\in(N,\,N+1)$, then
$M^{s,\,p}(\rn)=\ca F^{s}_{p,\,\fz}(\rn)$ with equivalent norms. Moreover,
$M^{s,\,p}(\rn)=F^{s}_{p,\,\fz}(\rn)$ when $s\in(N,\,N+1)$ and
$M^{N+1,\,p}(\rn)=F^{N+1}_{p,\,2}(\rn)$ with equivalent norms.
\end{cor}

\begin{rem}\label{r3.1} Notice that when $p\in(1,\,\fz)$, $h^p(\rn)=L^p(\rn)$,
and when $p\in(0,\,1]$, $h^p(\rn)\subsetneq L^p(\rn)$. Another way
to define the inhomogeneous Haj\l asz-Sobolev space denoted by $\wz
M^{s,\,p}(\rn)$ is to replace $\|f\|_{h^p(\rn)}$ by
$\|f\|_{L^p(\rn)}$ in the Definition \ref{d3.1}. Recall that it was
proved in \cite{ks} that $\wz M^{1,\,p}(\rn)= F^1_{p,\,2}(\rn)\cap
L^p(\rn)$ for $p\in(n/(n+1),\,\fz)$. An  argument similar to that
used in the proof of Theorem \ref{t3.2} can show that $\wz
M^{s,\,p}(\rn)=\ca^\ell_{0,\,m}F^s_{p,\,\fz}(\rn)\cap L^p(\rn)$ for
$p\in(n/(n+s),\,\fz)$ and $s\in(0,\,1]$. Similar results for
Corollary \ref{c3.x} also hold true. We omit the details.
\end{rem}

\section{Proofs of Theorems \ref{t1.3} and \ref{t1.4}}\label{s4}

\hskip\parindent The following spaces of test functions play a key
role in the theory of function spaces on RD-spaces; see \cite{hmy2}.
In what follows, for any $x,$ $y\in\cx$ and $r>0$, set $V(x,
y)\equiv\mu(B(x, d(x,y)))$ and $V_r(x)\equiv\mu(B(x, r))$. It is
easy to see that $V(x,y)\sim V(y,x)$ for all $x,\,y\in\cx$.

\begin{defn}\label{d4.1}
Let $x_1\in\cx$, $r\in(0, \fz)$, $\bz\in(0, 1]$ and $\gz\in(0,
\fz)$. A function $\vz$ on $\cx$ is said to be in the space
$\cg(x_1, r, \bz, \gz)$ if there exists a nonnegative constant
$C$ such that

(i) $|\vz(x)|\le C\frac{1}{V_r(x_1)+V(x_1, x)}
\lf(\frac{r}{r+d(x_1, x)}\r)^\gz$ for all $x\in\cx$;

(ii) $|\vz(x)-\vz(y)| \le C\lf(\frac{d(x, y)}{r+d(x_1, x)}\r)^\bz
\frac{1}{V_r(x_1)+V(x_1, x)} \lf(\frac{r}{r+d(x_1, x)}\r)^\gz$
for all $x$, $y\in\cx$ satisfying that $d(x, y)\le(r+d(x_1, x))/2$.

Moreover, for any $\vz\in\cg(x_1, r, \bz,\gz)$, its norm is defined
by $\|\vz\|_{\cg(x_1,\, r,\, \bz,\,\gz)} \equiv\inf\{C:\, (i) \mbox{
and } (ii) \mbox{ hold}\}$.
\end{defn}

Throughout the whole paper, we fix $x_1\in \cx $ and let $\cg(\bz,\gz)
\equiv\cg(x_1,1,\bz,\gz).$ Then
$\cg(\bz, \gz)$ is a Banach space.
We also let
$\ocg(\bz, \gz)=\lf\{f\in\cg(\bz, \gz):\,\int_\cx f(x)\,d\mu(x)=0\r\}.$
Denote by $(\cg(\bz, \gz))'$ and $(\ocg(\bz, \gz))'$
the dual spaces of $\cg(\bz, \gz)$ and $\ocg(\bz, \gz)$, respectively.
Obviously, $(\ocg(\bz, \gz))'=(\cg(\bz, \gz))'/\cc$.

For any given $\ez\in(0, 1]$, let $\cg_0^\ez(\bz, \gz)$ be the
completion of the set $\cg(\ez, \ez)$ in the space $\cg(\bz, \gz)$
when $\bz$, $\gz\in(0, \ez]$. Obviously, $\cg_0^\ez(\ez,
\ez)=\cg(\ez, \ez)$. If $\vz\in\cg_0^\ez(\bz, \gz)$, define
$\|\vz\|_{\cg_0^\ez(\bz, \gz)}\equiv\|\vz\|_{\cg(\bz, \gz)}$.
Obviously, $\cg_0^\ez(\bz, \gz)$ is a Banach space. The space
$\ocg_0^\ez(\bz,\gz)$ is defined to be the completion of the space
$\ocg(\ez,\ez)$ in $\ocg(\bz,\gz)$ when $\bz,\ \gz\in(0,\,\ez]$. Let
$(\cg_0^\ez(\bz, \gz))'$ and $(\cg_0^\ez(\bz, \gz))'$ be the dual
space of $\cg_0^\ez(\bz, \gz)$ and $\cg_0^\ez(\bz, \gz)$,
respectively. Also we have that
$(\ocg^\ez_0(\bz,\,\gz))'=(\cg^\ez_0(\bz,\,\gz))'/\cc$.

\begin{rem}\label{r4.1}
Because $(\ocg^\ez_0(\bz,\,\gz))'=(\cg^\ez_0(\bz,\,\gz))'/\cc$, if we
replace  $(\ocg^\ez_0(\bz,\,\gz))'$ with
$(\cg^\ez_0(\bz,\,\gz))'/\cc$ or $(\cg^\ez_0(\bz,\,\gz))'$
 in Definition \ref{d1.6}, then we obtain a new Triebel-Lizorkin space
which, modulo constants, is equivalent to the original
Triebel-Lizorkin space. So we can replace $(\ocg^\ez_0(\bz,\,\gz))'$
with $(\cg^\ez_0(\bz,\,\gz))'/\cc$ or $(\cg^\ez_0(\bz,\,\gz))'$ in
the Definition \ref{d1.6} if need be, in what follows.
\end{rem}

Now we recall the notion of approximations of the identity on
RD-spaces, which were first introduced in \cite{hmy2}.

\begin{defn}\label{d4.2}
Let $\ez_1\in(0, 1]$. A sequence $\{S_k\}_{k\in\zz}$ of bounded
linear integral operators on $L^2(\cx)$ is called an approximation
of the identity of order $\ez_1$ (for short, $\ez_1$-$\aoti$) with
bounded support, if there exist constants $C_3$, $C_4>0$ such that
for all $k\in\zz$ and all $x$, $x'$, $y$ and $y'\in\cx$, $S_k(x,
y)$, the integral kernel of $S_k$ is a measurable function from
$\cx\times\cx$ into $\cc$ satisfying
\begin{enumerate}
\item[(i)] $S_k(x, y)=0$ if $d(x, y)>C_42^{-k}$
and $|S_k(x, y)|\le C_3\frac{1}{V_{2^{-k}}(x)+V_{2^{-k}}(y)};$
\vspace{-0.3cm}
\item[(ii)] $|S_k(x, y)-S_k(x', y)|
\le C_3 2^{k \ez_1}[d(x, x')]^{\ez_1}
\frac{1}{V_{2^{-k}}(x)+V_{2^{-k}}(y)}$ for $d(x, x')\le\max\{C_4,
1\}2^{1-k}$; \vspace{-0.8cm}
\item[(iii)] Property (ii) holds with $x$ and $y$ interchanged;
\vspace{-0.3cm}
\item[(iv)] $|[S_k(x, y)-S_k(x, y')]-[S_k(x', y)-S_k(x', y')]|
\le C_3 2^{2k\ez_1}\frac{[d(x, x')]^{\ez_1}[d(y, y')]^{\ez_1}}
{V_{2^{-k}}(x)+V_{2^{-k}}(y)}$ for $d(x, x')\le\max\{C_4,
1\}2^{1-k}$ and $d(y, y')\le\max\{C_4, 1\}2^{1-k}$; \vspace{-0.3cm}
\item[(v)] $\int_\cx S_k(x, y)\, d\mu(y)
=1=\int_\cx S_k(x, y)\, d\mu(x)$.
\end{enumerate}
\end{defn}

It was proved in \cite[Theorem~2.6]{hmy2} that there always exists a
1-$\aoti$ with bounded support on an RD-space.

To prove Theorem \ref{t1.3}, we need a Sobolev embedding theorem,
which for $s=1$ is due to Haj\l asz \cite[Theorem 8.7]{h03}, and
for $s\in (0,1)$ can be proved by a slight modification of the
proof of \cite[Theorem 8.7]{h03}. We omit the details.

\begin{lem}\label{l4.1}
Let $s\in (0,1]$, $p\in (0,n/s)$ and
  $p^\ast=np/(n-sp)$.
Then there exists a positive constant $C$ such that for all
$u\in\dmsp$,  $g\in D(u)$ and all balls $B_0$ with radius $r_0$,
 $u\in L^{p^\ast}(B_0)$
and
\begin{equation*}
\inf_{c\in\rr}\lf(\bint_{B_0}|u-c|^{p^\ast}\,d\mu\r)^{1/p^\ast} \le
C r_0^s \lf(\dbint_{2B_0}g^p\,d\mu\r)^{1/p}.
\end{equation*}
\end{lem}

By Lemma \ref{l4.1}, we have the following version of Lemma
\ref{l2.2}.

\begin{lem}\label{l4.2}
Let $s\in(0,\,1]$, $p\in[n/(n+s),\,n/s)$
and $p_\ast=np/(n-sp)$.
Then for each $u\in \dot M^{s,\,p}(\cx)$,
there exists constant
$C$ such that
$u-C\in L^{p_\ast}(\cx)$
and $\|u-C\|_{L^{p_\ast}(\cx)}\le \wz C \|u\|_{\dot M^{s,\,p}(\cx)}$,
where  $\wz C$ is a positive
 constant independent of $u$ and $C$.
\end{lem}

With the aid of Lemmas \ref{l4.1} and \ref{l4.2}, we can prove
Theorem \ref{t1.3} by following the ideas used in the proof of Theorem
\ref{t1.1}. For reader's convenience, we sketch the argument.

\begin{proof}[Proof of Theorem \ref{t1.3}.]
We first prove that if $f\in \ca\dot F^s_{p,\,\fz}(\cx)$, then $f\in
\dot M^{s,\,p}(\cx)$ and $\|f\|_{\dot M^{s,\,p}(\cx)}\ls
\|f\|_{\ca\dot F^s_{p,\,\fz}(\cx)}$. Let $\{S_k\}_{k\in\zz}$ be a
$1$-$\aoti$ with bounded support. If $f$ is a locally integrable
function,  using $S_k(f)$ to replace $\vz_{2^{-k}}\ast f$ and
following the procedure as in the proof of Theorem \ref{t1.1}, we
know that $f\in \dot M^{s,\,p}(\cx)$ and
$$g(\cdot)=\sup_{k\in\zz}\sup_{\phi\in\ca_k(\cdot)}2^{ks}
|\langle f,\,\phi\rangle|\in L^p(\cx)$$ and $\|f\|_{\dot
M^{s,\,p}(\cx)}\ls \|f\|_{\ca\dot F^s_{p,\,\fz}(\cx)}$. If $f\in
\ca\dot F^s_{p,\,\fz}(\cx)$ is only known to be an element in
$(\cg(1,\,2))'$, we may also identify $f$ with a locally integrable
function $\wz f$ in $(\cg(1,\,2))'$ and $\|\wz f\|_{\dot
M^{s,\,p}(\cx)}\ls\|f\|_{\ca\dot F^s_{p,\,\fz}(\cx)}$ by using Lemma
\ref{l4.2} and an argument used in that of Theorem \ref{t1.1}. In
this sense, we have that $f\in \dot M^{s,\,p}(\cx)$ and $\|f\|_{\dot
M^{s,\,p}(\cx)}\ls\|f\|_{\ca\dot M^{s,\,p}(\cx)}$.

Conversely, let $f \in\dot M^{s,\,p}(\cx)$. Choose $g\in {\mathcal
D}(f)$ such that $\|g\|_{L^p(\cx)}\le 2\|f\|_{\dot M^{s,\,p}(\cx)}.$
Then for all $x\in\cx$, $k\in\zz$ and $\phi\in\ca_k(x)$,
 similarly to the proof of \eqref{e2.1} and using Lemma \ref{l4.1}, we have
that
\begin{eqnarray}\label{e4.5}
|\langle f,\,\phi\rangle|
&&\ls 2^{-ks}\lf[M\lf(g^{n/(n+s)}\r)(x)\r]^{(n+s)/n},\noz
\end{eqnarray}
which together with the $L^{p(n+s)/n}(\cx)$-boundedness of $M$ implies
that $\lf\|f\r\|_{\ca\dot F^s_{p,\,q}(\cx)} \ls\|g\|_\lp$ for all
$p\in(n/(n+s),\,\fz)$. Moreover,  without loss of generality, we may
assume that $M\lf(g^{n/(n+s)}\r)(x_1)<\infty$. Then for all $\psi\in\cg
(1,\,2)$, letting $\sz:=\int_\cx\psi(y)\,d\mu(y)$, by Lemma
\ref{l4.2} and an argument similar to the proof of \eqref{e2.1}, we
have that $f\in L^1_\loc(\cx)$ and
\begin{eqnarray*}
\lf|\int_\cx f(x)\psi(x)\,d\mu(x)\r|&&=\lf|\int_\cx f(x)
[\psi(x)-\sz S_0(x_1,\,x)]\,d\mu(x)\r|+|\sz||S_0(f)(x_1)|\\
&&\ls
 \|\psi\|_{L^1(\cx)}\dbint_{B(0,\,2C_4)}|f(z)|\,d\mu(z)+
 \|\psi\|_{\cg (1,\,2)}[M(g^{n/(n+s)})]^{1+s/n}(x_1)\\
&&
\ls C(f)\|\psi\|_{\cg (1,\,2)},
\end{eqnarray*}
which implies that $f\in(\cg (1,\,2))'$. Thus $f\in \ca\dot
F^s_{p,\,\fz}(\cx)$ and $\|f\|_{\ca\dot F^s_{p,\,\fz}(\cx)}
\ls\|f\|_{\dot M^{s,\,p}(\cx)},$ which completes the proof of
Theorem \ref{t1.3}.
\end{proof}

\begin{rem}\label{r4.2}
By the above proof, if we replace the space $\cg(1,\,2)$ of test functions by $\cg(\bz,\,\gz)$ with
 $\bz\in[s,\,1]$ and $\gz\in(s,\,\fz)$ in  Definition \ref{d1.5}, then Theorem \ref{t1.3} still holds true.
Thus the definition of the space $\ca\dot F^s_{p,\,\fz}(\cx)$ is independent of the choice of the
space of test functions $\cg(\bz,\,\gz)$  with
 $\bz\in[s,\,1]$ and $\gz\in(s,\,\fz)$.
\end{rem}

To prove Theorem \ref{t1.4}, we need the following homogeneous
Calder\'on reproducing formula established in \cite{hmy2}. We first
recall the following construction given by Christ in \cite{ch90},
which provides an analogue of the set of Euclidean dyadic cubes on
spaces of homogeneous type.

\begin{lem}\label{l4.3}
Let $\cx$ be a space of homogeneous type. Then there exists a
collection $\{Q^k_\az\subset\cx:\ k\in\zz,\ \az\in I_k\}$
of open subsets, where $I_k$ is some index set,
and constants $\dz\in (0,1)$ and $C_5,\ C_6>0$
such that
\begin{enumerate}
\vspace{-0.2cm}
\item[(i)] $\mu(\cx\setminus \cup_\az Q^k_\az)=0$ for each fixed $k$
and $Q^k_\az\cap Q^k_\bz=\emptyset$ if $\az\ne\bz$;
\vspace{-0.2cm}
\item[(ii)] for any $\az,\ \bz,\ k,\ \ell$ with $\ell\ge k,$ either $Q_\bz^\ell
\subset Q^k_\az$ or $Q^\ell_\bz\cap Q^k_\az=\emptyset$;
\vspace{-0.2cm}
\item[(iii)] for each $(k,\az)$ and each $\ell<k$,
there exists a unique $\bz$ such that $Q_\az^k\subset Q^\ell_\bz$;
\vspace{-0.2cm}
\item[(iv)] $\diam (Q_\az^k)\le C_5\dz^k$;
\vspace{-0.2cm}
\item[(v)] each $Q_\az^k$ contains some ball $B(z^k_\az, C_6\dz^k)$,
where $z^k_\az\in\cx$.
\end{enumerate}
\end{lem}

In fact, we can think of $Q^k_\az$ as being a {\it dyadic cube} with
diameter roughly $\dz^k$ and centered at $z^k_\az.$ In what follows,
to simplify our presentation, we always suppose that $\dz=1/2$;
otherwise, we need to replace $2^{-k}$ in the definition of
approximations to the identity by $\dz^k$ and some other changes are
also necessary; see \cite{hmy2} for more details.

In the following, for $k\in\zz$ and $\tau\in\ik$, we denote by
$Q_\tau^{k,\nu},$  $\nu=1,\ 2,\ \cdots, N(k,\tau),$ the set of all
cubes $Q_{\tau'}^{k+j}\subset Q_\tau^k,$ where $Q_\tau^k$ is the
dyadic cube as in Lemma \ref{l4.3} and $j$ is a fixed positive large
integer such that $2^{-j}C_5<1/3.$ Denote by $z_\tau^{k,\nu}$ the
``center" of $\qtn$ as in Lemma \ref{l4.3} and by $y_\tau^{k,\nu}$ a
point in $Q_\tau^{k,\nu}.$

\begin{lem}\label{l4.4}
Let $\ez\in (0,1)$ and $\{S_k\}_{k\in\zz}$ be a 1-$\aoti$ with
bounded support. For $k\in\zz$, set $D_k:=S_k-S_{k-1}$. Then, for any
fixed $j\in\nn$ large enough, there exists a family $\{\wz
D_k\}_{k\in\zz}$ of linear operators  such that for any fixed
$\ytn\in\qtn$ with $k\in\zz$, $\tau\in\ik$ and $\nu=1,\cdots,\nkt$,
$x\in\cx$, and all $f\in (\ocg_0^\ez(\bz,\gz) )'$ with $\bz,\
\gz\in(0,\,\ez)$,
$$f(x)=\dsum^\fz_{k=-\fz}\dsum_{\tau\in\ik}\dsum^\nkt_{\nu=1}
\mu(\qtn)\wz D_k(x,\ytn)D_k(f)(\ytn),$$
where the series converge in $ (\ocg_0^\ez(\bz,\gz) )'$.
Moreover, for any $\ez'\in (\ez,1)$, there exists a positive
constant $C$, depending on ${\ez'}$, such that the kernels,
denoted by $\wz D_k(x,y)$, of the operators $\wz D_k$ satisfy
\begin{enumerate}
    \item[(i)] for all $x,\ y\in\cx$,
    $|\wdk(x,y)|\le C\frac 1{V_{2^{-k}}(x)+V(x,y)}
    \lf[\frac {2^{-k}}{{2^{-k}}+d(x,y)}\r]^\ezp,$
    \item[(ii)] for all $x,\ x',\ y\in\cx$
with $d(x,x')\le ({2^{-k}}+d(x,y))/2$,
    \begin{eqnarray*}
    |\wdk(x,y)-\wdk(x',y)|&\le& C\lf[\frac {d(x,x')}{{2^{-k}}+d(x,y)}\r]^\ezp
    \frac 1{V_{2^{-k}}(x)+V(x,y)}\lf[\frac {2^{-k}}{{2^{-k}}+d(x,y)}\r]^\ezp,
\end{eqnarray*}
\item[(iii)] for all $k\in\zz$, $\int_\cx\wz D_k(x,y)\,d\mu(y) =0=\int_\cx
\wz D_k(x,y)\,d\mu(x).$
\end{enumerate}
\end{lem}

\begin{proof}[Proof of Theorem \ref{t1.4}.]
If $f\in \ca\dot F^s_{p,\,q}(\cx)$, by Remark \ref{r4.2} and the fact that
$ \ca\dot F^s_{p,\,q}(\cx)\subset \ca\dot F^s_{p,\,\fz}(\cx)$, we
know that $f\in (\cg(\bz,\,\gz))'$ with $\bz\in(s,\,1)$ and $\gz\in(s,\,\fz)$
and thus, $f\in \dot
F^s_{p,\,q}(\cx)$ and $\|f\|_{\dot F^s_{p,\,q}(\cx)}\ls
\|f\|_{\ca\dot F^s_{p,\,q}(\cx)}$. Conversely, assume that $f\in
\dot F^s_{p,\,q}(\cx)$. By Lemma \ref{l4.4}, for all $x\in\cx$,
$\ell\in\zz$\ and $\phi\in\ca_\ell(x)$,  we have
$$\langle f,\,\phi\rangle=\dsum^\fz_{k=-\fz}\dsum_{\tau\in\ik}\dsum^\nkt_{\nu=1}
\mu(\qtn)D_k(f)(\ytn)\int_\cx \wz D_k(z,\ytn)\phi(z)\,d\mu(z),$$ where we
fix $\ytn\in Q_\tau^{k,\,\nu}$ such that
\begin{equation}\label{e4.1}
|D_k(f)(\ytn)|\le 2\inf_{z\in Q_\tau^{k,\,\nu}} |D_k(f)(z)|.
\end{equation}
Recall that $\wz D_k$ depends on the choice of $\ytn$ and thus on
$f$, but they do have uniform estimates as in Lemma \ref{l4.4},
which is enough for us. In fact, by these estimates and
$\phi\in\ocg^\ez_0(\bz,\,\gz)$, we further know that for any fixed
$\bz'\in(s,\,\bz)$ and $\gz'\in(s,\,\gz)$ satisfying \eqref{1.6},
$$\lf|\int_\cx \wz D_k(z,\ytn)\phi(z)\,d\mu(z)\r|\ls 2^{-|k-\ell|\bz'}
\frac1{V_{2^{-(k\wedge\ell)}}(x)+V(x,\,\ytn)}
\lf[\frac{2^{-(k\wedge\ell)}}{2^{-(k\wedge\ell)}+d(x,\,\ytn)}\r]^{\gz'};$$
see \cite{hmy2} for a detailed proof. Thus,
choosing an $r\in(n/(n+[\bz'\wedge\gz']),\,\min\{p,\,q\})$, by
\eqref{e4.1}, we have
\begin{eqnarray*}
|\langle f,\,\phi\rangle|&&\ls
\dsum^\fz_{k=-\fz}2^{-|k-\ell|\bz'}\dsum_{\tau\in\ik}\dsum^\nkt_{\nu=1}
\frac{\mu(\qtn)|D_k(f)(\ytn)|}{V_{2^{-(k\wedge\ell)}}(x)+V(x,\,\ytn)}
\lf[\frac{2^{-(k\wedge\ell)}}{2^{-(k\wedge\ell)}+d(x,\,\ytn)}\r]^{\gz'}\\
&&\ls\dsum^\fz_{k=-\fz} 2^{-|k-\ell|\bz'}2^{[(k\wedge\ell)-k]n(1-1/r)}
\lf\{M\lf(\dsum_{\tau\in\ik}\dsum^\nkt_{\nu=1}
 |D_k(f)(\ytn)|^r\chi_{Q_\tau^{k,\,\nu}}\r)(x)\r\}^{1/r}\\
&&\ls\dsum^\fz_{k=-\fz} 2^{-|k-\ell|\bz'}2^{[(k\wedge\ell)-k]n(1-1/r)}
\lf[M\lf(|D_k(f)|^r\r)(x)\r]^{1/r}.
\end{eqnarray*}
This implies that
\begin{eqnarray*}
\|f\|_{\ca\dot F^s_{p,\,q}(\cx)}
&&\ls\lf\|\lf\{\sum_{\ell=-\fz}^\fz 2^{(\ell-k)sq}
\lf(
\dsum^\fz_{k=-\fz}2^{-|k-\ell|\bz'}2^{[(k\wedge\ell)-k]n(1-1/r)}\r.\r.\r.\\
&&\quad\quad\quad\quad\quad\quad\quad\quad\quad\quad
\lf.\lf.\lf.\times \lf[M\lf(2^{ksr}|D_k(f)|^r\r) \r]^{1/r}
 \r)^q
\r\}^{1/q}\r\|_{L^p(\cx)}.
\end{eqnarray*}
Applying the H\"older inequality when $q>1$ and the inequality that
$(\sum_{k}|a_k|)^q\le\sum_{k}|a_k|^q$ when $q\in(0,\,1]$ for all
$\{a_k\}_{k\in\zz}\subset\cc$, and using the vector-valued
inequality of the Hardy-Littlewood maximal operator (see
\cite{gly}), we then have
\begin{eqnarray*}
\|f\|_{\ca\dot F^s_{p,\,q}(\cx)}
&&\ls\lf\|\lf\{
\dsum^\fz_{k=-\fz}
\lf[M\lf(2^{ksr}|D_k(f)|^r\r)
\r]^{q/r}
\r\}^{1/q}\r\|_{L^p(\cx)}\ls\|f\|_{\dot F^s_{p,\,q}(\cx)}.
\end{eqnarray*}
This finishes the proof of Theorem \ref{t1.4}.
\end{proof}

\section{Inhomogeneous versions of Theorems \ref{t1.3} and
\ref{t1.4}}\label{s5}

\hskip\parindent We consider both cases $\mu(\cx)<\fz$
and $\mu(\cx)=\fz$ at the same time.
We next recall the notions of inhomogeneous Besov and Triebel-Lizorkin
spaces from \cite{hmy2}.

We call $\{S_k\}_{k\in\nn}$ to be an inhomogeneous approximation of
the identity of order $\ez$ with bounded support if their kernels
satisfy (i) through (v) of Definition \ref{d4.1}.

\begin{defn} \label{d5.1}
Let $\ez$, $s,\, p,\, q,\,\bz,\,\gz$ be as in Definition \ref{d1.6}.
Let $\{S_k\}_{k\in\nn}$ be an inhomogeneous approximation of the
identity of order $\ez$ with bounded support. For $k\in\nn$, set
$D_k\equiv S_k-S_{k-1}$. Let $\{\qto:\ \tau\in I_0,\ \nu=1,\cdots,
N(0,\tau)\}$ with a fixed large $j\in\nn$ be dyadic cubes as in
Section \ref{s4}. Let $s\in(0,\,\ez)$. The inhomogeneous
Triebel-Lizorkin space $\f$ is defined to be the set of all
$f\in\lf(\cg^\ez_0(\bz,\gz)\r)'$ that satisfy
\begin{eqnarray*}
\|f\|_\f&\equiv&\lf\{\dsum_{\tau\in I_0}\dsum^\nz_{\nu=1}
\mu(\qto)\lf[m_\qto(|S_0(f)|)\r]^p\r\}^{1/p}\\
&&+\lf\|\lf\{\dsum^\fz_{k=1} 2^{ksq}|D_k(f)|^q\r\}^{1/q}\r\|_\lp<\fz
\end{eqnarray*}
with the usual modification made when $q=\fz$.
\end{defn}

As shown in \cite{yz09}, the definition of $\f$
is independent of the choices of $\ez$, $\bz$, $\gz$
and the inhomogeneous approximation
of the identity.

\begin{defn}\label{d5.2}  Let $s\in(0,\,1]$,
$p\in(0,\,\fz)$ and $q\in(0,\,\fz].$ Let
$\ca:\equiv\{\ca_{k}(x)\}_{k\in\zz_+,\,x\in\cx}$ with
 $\ca_0(x)= \{\phi\in\cg(1,\,2),\,
\|\phi\|_{\cg(x,\,1,\,1,\,2)}\le1 \}$ and for $k\in\nn$,
$$\ca_k(x):= \{\phi\in\cg (1,\,2),\,
\|\phi\|_{\ocg(x,\,2^{-k},\,1,\,2)}\le1\}.$$
The inhomogeneous
grand Triebel-Lizorkin space $\ca\dot F^s_{p,\,q}(\cx)$
is defined to be the set of all
$f\in\lf(\cg(1,\,2)\r)'$ that satisfy
$$\|f\|_{\ca F^s_{p,\,q}(\cx)}\equiv\lf\|\lf\{\sum^\fz_{k=0}
2^{ksq}\sup_{\phi\in\ca_k(\cdot)}
|\langle f,\,\phi\rangle|^q\r\}^{1/q}\r\|_\lp<\fz$$
with the usual modification made when $q=\fz$.
\end{defn}

Then we have the following result.

\begin{thm}\label{t5.1}
Let all the assumptions be as in Definition \ref{d5.1}.
Then $F^s_{p,\,q}(\cx)=\ca F^s_{p,\,q}(\cx)$ with equivalent norms.
\end{thm}

The proof of Theorem \ref{t5.1} is similar to that of Theorem
\ref{t1.4}. We point out that instead of the homogeneous Calder\'on
reproducing formula, we need the inhomogeneous one established in
\cite{hmy2}. We omit the details.

Moreover, we define the inhomogeneous Haj\l asz-Sobolev spaces as
follows. We also have an inhomogeneous version of Theorem
\ref{t1.3}, which can be proved by using the ideas appearing in the
proofs of Theorem \ref{t3.2} and Theorem \ref{t1.3}.

\begin{defn}\label{d5.3} Let $p\in (0,\fz)$ and $s\in (0,1]$. The
inhomogeneous fractional Haj\l asz-Sobolev space $\msp$ is defined
to be the set of all measurable functions $f\in L^p_\loc(\cx)$ that satisfy
both  $f\in h^p(\cx)=F^0_{p,\,2}(\cx)$ and $f\in\dmsp$; moreover,
define $\|f\|_\msp\equiv \|f\|_{h^p(\cx)}+\inf_{g\in
\cd(f)}\|f\|_\dmsp,$
\end{defn}

\begin{thm}\label{t5.2}
Let $s\in (0,1]$ and $p\in(n/(n+s), \fz)$. Then
$\msp=\ca F^s_{p,\,\fz}(\cx)$ with equivalent norms.
\end{thm}

\bigskip

\noindent {\sc Pekka Koskela}: University of Jyv\"askyl\"a,
Department of Mathematics and Statistics, P.O. Box 35 (MaD),
Fin-40014 University of Jyv\"askyl\"a, Finland

\medskip

\noindent{\it E-mail:} \texttt{pkoskela@maths.jyu.fi}

\bigskip

\noindent {\sc Dachun Yang}: School of Mathematical Sciences,
Beijing Normal University,
Laboratory of Mathematics and Complex Systems,
Ministry of Education,
Beijing 100875, People's Republic of China

\medskip

\noindent{\it E-mail:} \texttt{dcyang@bnu.edu.cn}

\bigskip

\noindent {\sc Yuan Zhou}: School of Mathematical Sciences, Beijing
Normal University, Laboratory of Mathematics and Complex Systems,
Ministry of Education, Beijing 100875, People's Republic of China

and

\noindent University of Jyv\"askyl\"a, Department of Mathematics and
Statistics, P.O. Box 35 (MaD), Fin-40014 University of
Jyv\"askyl\"a, Finland

\medskip

\noindent{\it E-mail:} \texttt{yuanzhou@mail.bnu.edu.cn}

\end{document}